\documentclass[11pt]{amsart}
\usepackage{amsmath,amssymb,mathrsfs,amscd}
\usepackage{graphicx}
\usepackage[all]{xy}

\title{A new series of compact minitwistor spaces and Moishezon twistor spaces over them}

\author{Nobuhiro Honda}
\thanks
{$^{\dag}$This work was partially supported by
Research Fellowships of the 
Japan Society for the Promotion
of Science for Young Scientists.
}
\date{}

\newcommand{\ol}{\overline}

\newcommand{\lra}{\longrightarrow}

\newcommand{\set}{\,|\,}
\newcommand{\proofend}{\hfill$\square$}

\newtheorem{prop}{Proposition}[section]
\newtheorem{lemma}[prop]{Lemma}

\newtheorem{thm}[prop]{Theorem}

\newtheorem{cor}[prop]{Corollary}

\newtheorem{definition}[prop]{Definition}

\begin{document}

\begin{abstract}

In recent papers \cite{Hon07-3, Hon07-4}, we gave explicit description of some new Moishezon twistor spaces.
In this paper, developing the method in the papers much further, we explicitly give projective models of a number of new Moishezon twistor spaces, as conic bundles over some rational surfaces (called minitwistor spaces). 
These include the twistor spaces studied in the papers as  very special cases.

Our source of the result is a series of self-dual metrics with torus action constructed by D.\,Joyce \cite{J95}.
Actually, for arbitrary Joyce metrics and $U(1)$-subgroups of the torus which fixes a torus-invariant 2-sphere,  we first determine the associated minitwistor spaces in explicit forms. 
Next by analyzing  the meromorphic maps from the twistor spaces to the minitwistor spaces, we realize projective models of the twistor spaces of all Joyce metrics, as conic bundles over the minitwistor spaces.
Then we prove that for any one of these minitwistor spaces, there exist Moishezon twistor spaces  
{\em with only $\mathbf C^*$-action} whose quotient space is the given minitwistor space.
This result generates numerous Moishezon twistor spaces which cannot be found in the literature
(including the author's papers), in quite explicit form.
\end{abstract}

\maketitle

\bigskip\noindent
\section{ Introduction}
By the fundamental theorem of the twistor theory initiated by R. Penrose, there is a natural one-to-one correspondence between self-dual conformal classes on a 4-manifold, and twistor spaces over the 4-manifold.
 N.\,J.\,Hitchin \cite{Hi81} showed that there is a different but similar one-to-one correspondence between Einstein-Weyl structures on  3-manifolds, and certain complex surfaces called {\em minitwistor spaces}.
P.\,E.\,Jones--K.\,P.\,Tod \cite{JT85} pointed out that, when the self-dual conformal class admits an action of a 1-dimensional Lie group, the former correspondence   naturally induces  the latter one.
Namely, an Einstein-Weyl structure is naturally induced on the orbit space of the action, and a minitwistor space can be obtained as a quotient space of the twistor space with respect to a lifted and complexified action.
In particular, if one finds a minitwistor space (or an Einstein-Weyl space), it may be possible to construct a twistor space (or a self-dual conformal class) with 1-dimensional symmetries over the minitwistor space (or the Einstein-Weyl space).
This idea is successfully realized by C.\,LeBrun \cite{LB91}, who constructed  explicit examples of self-dual metrics and their twistor spaces on compact 4-manifolds, built over a basic example of Einstein-Weyl space (the upper-half space $\mathscr H^3$) and its minitwistor space (the product surface $\mathbf{CP}^1\times\mathbf{CP}^1$) .

In \cite{Hon07-3, Hon07-4}, pursuing this idea, the author  constructed some explicit examples of compact but singular minitwistor spaces and then constructed  Moishezon twistor spaces whose projective models have a natural structure of conic bundles over the minitwistor spaces.
In this paper, developing this method much further, we systematically provide a number of minitwistor spaces and Moishezon twistor spaces  with $\mathbf C^*$-action over them.
These include the twistor spaces in \cite{Hon07-3, Hon07-4} as very special cases.  (However  much more investigations are required  to reach the detailed results about the structure and the construction given in \cite{Hon07-3} and \cite{Hon07-4}).

A resource of the present construction is  a series of self-dual metrics with 2-dimensional torus action constructed by D.\,Joyce \cite{J95}, and their associated twistor spaces investigated by A.\,Fujiki \cite{F00}.
Namely, in Section 2, we start from an arbitrary Joyce metric on $n\mathbf{CP}^2$ (the connected sum of $n$ copies of complex projective planes) and choose an arbitrary $U(1)$-subgroup of the torus that fixes a 2-sphere in $n\mathbf{CP}^2$.
In other words, we take an arbitrary $\mathbf C^*$-subgroup of $\mathbf C^*\times\mathbf C^*$ that fixes a conjugate pair of  smooth rational curves on the twistor space.
For these subgroups, {\em we can explicitly give a linear system on the twistor space whose associated meromorphic map can be regarded as a quotient map of the $\mathbf C^*$-action}
(Proposition \ref{prop-mt3}).
We can also determine the structure of the image surfaces (namely the minitwistor spaces) in completely explicit form (Propositions \ref{prop-mt5}, \ref{prop-rf1}, \ref{prop-sing_T}).
We emphasizes that these minitwistor spaces are given as  complex surfaces embedded in projective spaces, rather than abstract complex surfaces.

By our construction, there is a natural meromorphic map from the  minitwistor spaces to  rational normal curves whose fibers are conics.
The location of the reducible fibers of this rational conic bundle map depends on the continuous parameters involved in the Joyce metrics (Proposition \ref{prop-rf1}).
Consequently, {\em  our minitwistor spaces constitute  moduli spaces, in accord with a variation of Joyce's parameters}.
But even if we neglect the difference raised by these deformations, {\em our method yields an infinite (but countable) number of minitwistor spaces}, corresponding to the choice of $n$,  torus actions on $n\mathbf{CP}^2$, and  $U(1)$-subgroups of the torus.
We remark that most of these minitwistor spaces are new, as far as the author knows.
Another interesting feature of these minitwistor spaces is that, except the simplest example (which is a smooth quadratic surface in $\mathbf{CP}^3$, the minitwistor space of the hyperbolic space), they always have some mild singularities (Proposition \ref{prop-sing_T}).

Since the linear systems on the twistor spaces inducing the quotient maps to the minitwistor spaces always have  base points, the maps have indeterminacy loci.
The structure of this base locus is so complicated in general that it looks difficult to give an explicit sequence of blowing-ups which removes the indeterminacy  completely.
But it is possible to give partial eliminations explicitly, which are enough for the purpose of obtaining  resolutions of the singularities of the minitwistor spaces, as well as  explicit $\mathbf{CP}^2$-bundles over the resolved minitwistor spaces (\S 3.1).
Then {\em projective models of the twistor spaces of (arbitrary) Joyce metrics can be realized as  conic bundles in these $\mathbf{CP}^2$-bundles} (\S 3.2).
These can be regarded as a generalization of  LeBrun's realization of (projective models of) his twistor spaces as conic bundles over $\mathbf{CP}^1\times\mathbf{CP}^1$,
because the last $\mathbf{CP}^1\times\mathbf{CP}^1$ is the most simple example of the minitwistor spaces we obtained.

So far we have discussed  the structure of the twistor spaces of Joyce metrics themselves. 
In order to obtain new twistor spaces, we consider  $\mathbf C^*$-equivariant deformations of these twistor spaces, where the $\mathbf C^*$-action is the one we have been considering.
For some of these $\mathbf C^*$-actions, any equivariant deformations become $\mathbf C^*\times\mathbf C^*$-equivariant and in that case we cannot obtain new twistor spaces.
However, we can give a sufficient condition for the $\mathbf C^*$-subgroup to have an equivariant deformation with the properties that (i)  not the whole of the $\mathbf C^*\times\mathbf C^*$-action survives but only the $\mathbf C^*$-action survive under the deformation, (ii) the structure of the quotient space by the $\mathbf C^*$-action remains unchanged under the deformation.
By the property (i), the deformed twistor spaces are not the twistor spaces of Joyce metrics,
but they are yet Moishezon by  the property (ii).
The sufficient condition is expressed in terms of some integer easily determined by the $U(1)$-action (Theorem \ref{thm-1}).
By using this criterion, we can further show that {\em for any prescribed minitwistor spaces obtained in Section 2, there always exist  Moishezon twistor spaces  with $\mathbf C^*$-action whose minitwistor space is the given one and whose automorphism group is just $\mathbf C^*$} (Theorem \ref{thm-2}).
We also determine discriminant locus of the quotient map from these twistor spaces to the minitwistor spaces (Proposition \ref{prop-d5}).
Hence we obtain the structure of these deformed twistor spaces in quite detailed form.

These results  produce a huge number of Moishezon twistor spaces with $\mathbf C^*$-action in a systematic way.
To justify, let $\delta(n)$ be the number of equivalent classes of $U(1)$-actions on $n\mathbf{CP}^2$ that can be obtained from  effective torus-actions on $n\mathbf{CP}^2$ by taking a $U(1)$-subgroup that fixes one of the torus-invariant spheres.
(Here we are considering all effective torus actions on $n\mathbf{CP}^2$. We have $\delta(1)=1,\,\delta(2)=2,\,\delta(3)=3,\,\delta(4)=7,\,\delta(5)=15$, etc.)
Then if $n\ge 3$ {\em there exist at least $\delta(n-1)$ effective $U(1)$-actions on $n\mathbf{CP}^2$  for which there exist non-Joyce  self-dual metrics invariant under the $U(1)$-action  
 whose twistor spaces are Moishezon} (Corollary \ref{cor-1}).
We note that the number of previously known $U(1)$-actions satisfying these properties was {\em only a few} for each $n$, and that this result is not known even if we do not require for the twistor spaces to be Moishezon.
 We also remark that thanks to Proposition \ref{prop-d5} mentioned above,
this is not just the existence theorem, but we have their detailed structure of the twistor spaces in  hand.

In Section 5.2--5.4,
in order to obtain explicit new examples, we give various $U(1)$-actions on $n\mathbf{CP}^2$.
In Section 5.2, we consider all such $U(1)$-actions 
 that can be obtained from LeBrun's torus-action. 
From these $U(1)$-actions we can already produce a considerable number of new Moishezon spaces with $\mathbf C^*$-action.
In Section 5.3, we first display all such $U(1)$-actions whose non-trivial isotropy subgroup 
is $\{\pm 1\}\subset U(1)$.
These $U(1)$-actions might be placed to a position next to the semi-free action.
We show that most of them have new Moishezon twistor spaces.
In Section 5.4, we consider a particular series of $U(1)$-actions on $n\mathbf{CP}^2$ which can be placed to another extreme
(in comparison with the $U(1)$-actions in \S 5.3).
The minitwistor spaces associated to these $U(1)$-actions are most singular ones among all minitwistor spaces we construct in this paper.
From these $U(1)$-actions we can also produce a lot of new Moishezon twistor spaces.
Meanwhile, we also arrange the twistor spaces obtained in \cite{Hon07-3} and \cite{Hon07-4} into the present framework.

 
 Apart from these explicit examples of new Moishezon twistor spaces,
an interesting observation concerning our minitwistor spaces is that,  
some hyperplane sections appear as discriminant curves of the conic bundle structure
 (Proposition \ref{prop-d5}).
 This is also true for the image of general twistor lines.
 Namely, so called the `minitwistor lines' are hyperplane sections of the minitwistor spaces.
This seems to indicate that  our realization of minitwistor spaces and twistor spaces is a natural one, and also that the minitwistor spaces should be considered not as abstract complex surfaces but as embedded surfaces in projective spaces, as we do in Section 2.

As the twistor spaces of Joyce metrics are Moishezon manifolds, it is possible to consider `canonical' quotient spaces under the $\mathbf C^*$-actions.
In Appendix we show that our minitwistor spaces are actually isomorphic to these canonical quotient spaces (Proposition \ref{prop-cq}).
This gives an intrinsic characterization of our minitwistor spaces.

The author expresses his gratitude to Akira Fujiki for explaining him results on the existence of  `canonical' quotient spaces by meromorphic actions on  manifolds in $\mathscr C$ in general, and asking him 
whether the author's minitwistor spaces are actually identical to the canonical quotient spaces,
which results in the content in Appendix.

\vspace{2mm}
{\bf Notations and Conventions.}
As in \cite{Hon07-3, Hon07-4}, to save  notations we adapt the following convention.
If $\mu:X\to Y$ is a bimeromorphic morphism of  complex variety and $W$ is a complex subspace in $X$, we write $W$ for the image $\mu(W)$  if the restriction $\mu|_W$ is still bimeromorphic.
Similarly, we use the same notation for a complex subspace in $Y$ and its strict transform into $X$.
If $D$ is a divisor on a variety $X$, the dimension of a complete linear system $|D|$ always means $\dim H^0(X,[D])-1$.
The base locus is denoted by ${\rm Bs}\,|D|$.
If a Lie group $G$ acts on $X$ by means of biholomorphisms and $D$ is $G$-invariant,
$G$ naturally acts on the vector space $H^0(X,[D])$.
Then $H^0(X,[D])^{G}$ means the subspace of all $G$-invariant sections.
Further $|D|^G$ means its associated linear system.
$X^G$ means the set of $G$-fixed points.
If $Z$ is a twistor space, $F$ always denotes the canonical square root of the anticanonical line bundle of $Z$. 
The degree of a divisor on $Z$ means its intersection number with twistor lines.
Throughout Sections 2--5, $K$ and $G$ mean the 2-dimensional Lie group $U(1)\times U(1)$ and  its complexification $\mathbf C^*\times\mathbf C^*$ respectively.

\section{Minitwistor spaces associated to the twistor spaces of Joyce metrics}
 
 We first recall basic properties of Joyce metrics and their twistor spaces.
For each $n\ge 1$, Joyce metrics on $n\mathbf{CP}^2$ are uniquely determined by the following two data.
One is a diffeomorphism type of an effective $K$-action on  $n\mathbf{CP}^2$, where the $K$-action becomes  the conformal isometry group of Joyce metrics.
Here, we say that two $K$-actions on  $n\mathbf{CP}^2$  have the same diffeomorphism type if there exists a $K$-equivariant diffeomorphism of $n\mathbf{CP}^2$.
For each $n$, the number of diffeomorphic types of effective $K$-actions on $n\mathbf{CP}^2$ is finite.
The other data is a set of different $(n+2)$ real numbers which actually determines a Joyce metric.
For fixed diffeomorphism type of $K$-action on $n\mathbf{CP}^2$, 
two Joyce metrics on $n\mathbf{CP}^2$  determined by   $\{\lambda_1<\lambda_2<\cdots<\lambda_{n+2}\}$ and $\{\lambda'_1<\lambda'_2<\cdots<\lambda'_{n+2}\}$ are mutually conformally isometric iff there exist real numbers $a,b,c,d$ with $|ad-bc|=\pm 1$ such that $\lambda_j'=(a\lambda_i+b)/(c\lambda_i+d)$ holds for any $1\le i\le n+2$ and some $1\le j\le n+2$.
Following Fujiki \cite[\S 3]{F00}, we call these $(n+2)$ numbers as the {\em conformal invariant} of Joyce metrics.

By a theorem of A.\,Fujiki \cite[Theorem 1.1]{F00}, if a self-dual conformal class on $n\mathbf{CP}^2$ admits an effective $K$-action of conformal isometries, then it coincides with the conformal class of some Joyce metric.
In the course of the proof, he showed the following result on the structure of the twistor spaces of Joyce metrics.
 
 \begin{prop}\label{prop-01}\cite{F00}
Let $Z$ be the twistor space of a Joyce metric on $n\mathbf{CP}^2$.
Then  the following holds.
(i) $\dim |F|^G=1$. 
Moreover, general members of the pencil are isomorphic to a smooth toric surface satisfying  $c_1^2=8-2n$.
(ii) Let $S$ be a general member of the pencil 
and  $C$ the  unique $G$-invariant anticanonical curve on $S$, which is a cycle of $2(n+2)$ rational curves.
Then {\rm{Bs}}\,$|F|^G=C$.
(iii) If we write $C=\sum_{i=1}^{2n+4}C_i$ in such a way that $C_i$ and $C_{i+1}$ intersect, the real structure on $Z$ exchanges $C_i$ and $C_{i+n+2}$, where the subscripts are counted modulo $(2n+4)$.
(iv) The pencil $|F|^G$ has precisely $(n+2)$ singular members and all of them consist of two smooth irreducible toric surfaces that are conjugations of each other.
(v) If we write $S_i=S_i^++S_i^-$ ($1\le i\le n+2$) for the reducible members, $S^+_i$ and $S^-_i$ divide $C$ into `halves' in the sense that both $S_i^+\cap C$ and $S_i^-\cap C$ are connected (cf.\! (iii)). 
Moreover,   if we write $L_i:=S^+_i\cap S_i^-$, then $\{L_i\set1\le i\le n+2\}$ coincides with the set of  $G$-invariant twistor lines on $Z$. 
\end{prop}

In the proposition, if a Joyce metric is not  a LeBrun metric (with $K$-action), we can replace $|F|^G$ by just $|F|$.
Note also that the toric surface $S$ is uniquely determined by the diffeomorphism type of $K$-action on $n\mathbf{CP}^2$. 
The conformal invariant corresponds to the place of reducible fibers of $|F|^G$.

 Let $S$ and $C$ be as in the proposition. By (iii), we can write the cycle $C$ as
\begin{equation}\label{cycle1}
C=\sum_{1\le i\le n+2} C_i+\sum_{1\le i\le n+2} \ol{C}_i.
\end{equation}
Here we are giving a numbering satisfying $C_i C_{i+1}=\ol{C}_i\ol{C}_{i+1}=1$ for $1\le i\le n+1$. Of course, this numbering is determined only up to cyclic permutations (applied to the whole of $C$) and reversing the orientation.
We make a distinction between $S_i^+$ and $S_i^-$ by declaring that $S_i^-$ contains the component $C_1$.

In the following we choose and fix {\em any} one of the irreducible components of $C$, and adapt a numbering for the irreducible components such that the chosen component is $C_1$.
Let $G_1\,(\simeq\mathbf C^*)$ be the isotropy subgroup of $C_1$.
Then the conjugate component $\ol{C}_1$ is also fixed by $G_1$.
Since the two curves $C_1$ and $\ol{C}_1$ are fixed, there is a natural holomorphic quotient map
\begin{equation}\label{quot_S}
S\,\lra\,\mathbf{CP}^1
\end{equation}
 by the $G_1$-action on $S$, 
which has $C_1$ and $\ol{C}_1$ as  distinguished sections.
On the base $\mathbf{CP}^1$, the quotient group $G/G_1\simeq\mathbf C^*$ naturally acts non-trivially.
This action has two fixed points which are conjugate of each other.
Let $f$ and $\ol{f}$ be the fibers (as schemes) of the quotient map \eqref{quot_S} over these two points.
These are the only reducible fibers of the quotient map \eqref{quot_S}.
As curves on $S$,  $f$ and  $\ol{f}$  are linearly equivalent.
These are non-reduced curves in general.
Actually,  writing these as
\begin{equation}\label{f_1}
f=\sum_{2\le i\le n+2}k_iC_i\,\,\text{ and }\,\,{\ol f}=\sum_{2\le i\le n+2}k_i\ol{C}_i,
\end{equation}
we have $k_2=k_{n+2}=1$ (since $f\cdot C_1=f\cdot \ol{C}_1=1$) but $k_i\ge 1$ for $2< i< n+2$.
In this way, for the prescribed component $C_1$ or the subgroup $G_1$, we obtain the sequence of $(n+1)$ positive integers $(k_2,k_3,\cdots,k_{n+2})$, where $k_2=k_{n+2}=1$. 
From this, we derive a positive integer $m$ through the following procedure:

\begin{enumerate}
\item[(i)]
Let $k$ be the biggest number among $\{k_2,k_3,\cdots,k_{n+2}\}$, and $i$ the smallest index satisfying $k_i=k$.
Next let $j$ be the biggest index satisfying $k_i=k_{i+1}=\cdots=k_j=k$.
(So $j\ge i$ holds. If $k_i=k$ is an `isolated' maximum, $j=i$ holds.)
\item[(ii)]
For each index $l$ satisfying $i\le l\le j$, we define $k'_l=k_l-1$.
Then consider the new sequence of $(n+1)$ integers obtained from the original sequence  $(k_2,k_3,\cdots,k_{n+2})$ by replacing $k_l$ by $k'_l$ for $i\le l\le j$.
\item[(iii)] \label{item-3}
If all the entries of the resulting sequence is zero, then the procedure ends.
If not, return to (i) and apply the same procedure.
\end{enumerate}
In the following we refer this procedure so frequently that we call as {\em Procedure (A)}.

\begin{definition}
\label{def-m}{\rm{
For the toric surface $S\in |F|^G$ and the subgroup $G_1\subset G$, let $(k_1,k_2,\cdots,k_{n+2})$ be the sequence representing the reducible fiber as above. Then we define $m$ to be the number of times we require Procedure (A) until it end.}}
\end{definition}
For instance, if $n=4$ and the sequence is $(1,2,1,2,1)$, then we have $m=3$.
Since each entries $k_i$ decreases by at most
one for each application of Procedure (A), we have 
 the relations
\begin{equation}\label{ineq-1}
m\ge k_i
\end{equation}
for all $2\le i\le n+2$.
In this way for the toric surface $S$ (in the twistor spaces of Joyce metrics) and each $\mathbf C^*$-subgroup $G_1$ fixing a component of the anticanonical cycle, we have associated a positive integer $m$.
This integer will be a basic invariant in our construction of  minitwistor spaces.
Note that once the pair $(S,G_1)$ is  given, it is very easy to determine the sequence $(k_2,k_3,\cdots,k_{n+2})$ and  to compute the integer $m$.

Among all twistor spaces of Joyce metrics, LeBrun twistor spaces  can be readily characterized in terms of the invariant $m$ as follows.

\begin{prop}\label{prop-char-LB}
Let $Z$ be the twistor space of a Joyce metric.
Then $Z$ is a LeBrun twistor space with $G$-action iff there exists a $\mathbf C^*$-subgroup $G_1$ of $G$ such that the integer $m$ for the pair  $(S,G_1)$ is one, where $S\in |F|^G$ is a smooth member.
\end{prop}

Since we do not use this result, we only mention that,
by a theorem of LeBrun \cite{LB93} on a characterization of his metrics in terms of $U(1)$ or $\mathbf C^*$-action, a key of the proof is the fact that $m=1$ if and only if $G_1$ is a subgroup acting semi-freely.

%


The following proposition gives a linear system on the toric surface $S$ which induces the $G_1$-quotient map \eqref{quot_S}.
\begin{prop}\label{prop-mt2}
Let $(S,G_1)$  be a pair as above and $m$ the positive integer determined in Definition \ref{def-m}. Then the linear subsystem (of the system $|mK_S^{-1}|$) generated by the 3 curves
\begin{equation}\label{gen-1}
mC,\,\,mC+f-{\ol f},\,\,mC-f+{\ol f}
\end{equation}
satisfies the following.
(Note that the latter 2 curves are  effective by the inequality \eqref{ineq-1}.)
(i)
Its movable part is free
(ii) The image of the associated holomorphic map is a conic in $\mathbf{CP}^2$.
(iii)  The holomorphic map coincides with the $G_1$-quotient map \eqref{quot_S}.
\end{prop}

\noindent
{\bf Proof.}
This is entirely elementary.
Actually, removing the maximal common curve $mC_1+m\ol C_1+\sum_{i=2}^{n+2}(m-k_i)(C_i+\ol{C}_i)$ from  \eqref{gen-1}, we obtain the 3 curves 
 \begin{equation}
\sum_{2\le i\le n+2}k_iC_i+\sum_{2\le i\le n+2}k_i\ol{C}_i,\,\,
\sum_{2\le i\le n+2}2k_iC_i,\,\,\sum_{2\le i\le n+2}2k_i\ol{C}_i.
\end{equation}
respectively.
These are exactly $f+{\ol f},2f$ and $2{\ol f}$ respectively.
From this all the claims (i)--(iii) follow immediately.
\proofend

\vspace{2mm}
With this proposition, for the twistor spaces of any Joyce metrics,  we are going to obtain generators of a subsystem of   $|mF|$ whose associated map gives a (meromorphic) quotient map for the twistor space with respect to the prescribed subgroup $G_1$.
($m$ will be the integer in Definition \ref{def-m}.)
This plays a key role for obtaining the minitwistor spaces.

\begin{prop}\label{prop-mt3}
Let $Z$ be the twistor space of a Joyce metric on $n\mathbf{CP}^2$, and $S$ a smooth $G$-invariant member of the pencil $|F|^G$.
Let $G_1\subset G$ be the isotropy subgroup of a component $C_1$ of the anticanonical cycle $C$ in $S$.
Let $m$ be the positive integer for  $(S,G_1)$ as in Definition \ref{def-m}.
Then we can explicitly find a $G$-invariant divisor $Y\in |mF|$  satisfying the following.
(i) All irreducible components of $Y$ are of degree one.
(ii) $Y|_S=mC+f-{\ol f}$ holds, where $Y|_S$ means the restriction onto $S$.
(iii) As irreducible components, $Y$ contains  $S_1^+$ and $S_{n+2}^-$ by multiplicity exactly one for each, and does not contain $S_1^-$ nor $S_{n+2}^+$.
In particular, $Y$ is not real.
(iv)  $Y$ does not contain  $S_j^+$ nor $S_j^-$ at the same time for any $1\le j\le n+2$.
\end{prop}

\noindent
{\bf Proof.} 
We give explicit way how to obtain the divisor $Y$.
First recall that we have
\begin{equation}\label{rest-1}
mC+f-{\ol f}=mC_1+\sum_{2\le i\le n+2}(m+k_i)C_i+m\ol{C}_1+\sum_{2\le i\le n+2}(m-k_i)\ol{C}_i,
\end{equation}
where we are still representing the fiber $f$ as in \eqref{f_1}.
For each integer $l$ satisfying  $1\le l\le m$ we explicitly take two degree one divisors as follows.
For each such an $l$, let $i_l$ and $j_l$ be the two indices $i$ and $j$ respectively obtained by the step (i) of $l$-th application of Procedure (A).
So we have $2\le i_l\le j_l\le n+2$.
We choose two degree one divisors $S_{i_l-1}^+$ and $S_{j_l}^-$.
Then  the restrictions of these onto $S$ are `connected halves of the cycle $C$' of the following form:
\begin{equation}\label{rest-37}
S_{i_l-1}^+|_S=C_{i_l}+C_{i_l+1}+\cdots+\ol{C}_{i_l-1},\,\,\,
S_{j_l}^-|_S=C_{j_l}+C_{j_l-1}+\cdots+\ol{C}_{j_l+1}.
\end{equation}
Take the sum of all these and define 
\begin{equation}\label{Y}
Y=\sum_{1\le l\le m}\left(S_{i_l-1}^++S_{j_l}^-\right).
\end{equation}
The degree of $Y$ is clearly $2m$.
(Since $i_l=i_{l'}$ or $j_l=j_{l'}$ can occur for $l\neq l'$, $Y$ is non-reduced in general.)
Then we have
\begin{align}
Y|_S&=\sum_{1\le l\le m}S_{i_l-1}^+|_S
\,\,+\sum_{1\le l\le m}S_{j_l}^-|_S\\
&=\sum_{1\le l\le m}(C_{i_l}+C_{i_l+1}+\cdots+\ol{C}_{i_l-1})
+\sum_{1\le l\le m}(C_{j_l}+C_{j_l-1}+\cdots+\ol{C}_{j_l+1}).\label{rest-28}
\end{align}
We now claim that for each $2\le i\le n+2$, in \eqref{rest-28}, $C_i$ and $\ol C_i$ are contained by multiplicities $(m+k_i)$ and $(m-k_i)$ respectively.
Actually, in the $l$-th application of Procedure (A), 
the component $C_i$ is chosen twice (i.e. both $S_{i_l-1}^+$ and $S_{j_l}^-$ contain $C_i$) iff $i_l\le i\le j_l$, and  otherwise chosen precisely once
(i.e. only one of $S_{i_l-1}^+$ and $S_{j_l}^-$ contains $C_i$).
On the other hand, in the $l$-th application,
the component $\ol C_i$ is not chosen at all iff $C_i$ is chosen twice and otherwise  precisely once (by \eqref{rest-37}).
Since the number of $l$ satisfying $i_l\le i\le j_l$ has to be exactly $k_i$,
the claim follows.
Then by \eqref{rest-1}, this means $Y|_S=mC+f-{\ol f}$.
Further, since $\ol{C}=C$, $\ol{Y}|_S=mC-f+{\ol f}$ holds.
Hence by Proposition \ref{prop-mt2}, $Y|_S$ and $\ol{Y}|_S$ belong to the system $|mK_S^{-1}|$.
Since the restriction map $H^2(Z,\mathbf Z)\to H^2(S, \mathbf Z)$ is always injective, this means that $Y$ and $\ol{Y}$ belong to the system $|mF|$.
By the choice, (i) and (ii) are clear.

Next we show (iii).
To see this, we notice that after applying Procedure (A) $(m-1)$ times, all entries of the resulting sequence become $1$.
Thus for the final application of (A), we have $i=2$ and $j=n+2$.
Namely in this final application, $S_1^+$ and $S_{n+2}^-$ are chosen.
Further, since $k_2=k_{n+2}=1$ as is already remarked,
these 2 divisors are not chosen during the preceding $(m-1)$ applications.
On the other hand, by our choice the degree one divisors for $Y$, it is obvious that $S_1^-$ and $S_{n+2}^+$ are never chosen.
Hence $Y$ contains $S_1^+$ and $S_{n+2}^-$ by multiplicity one, and does not contain $S_1^-$ nor $S_{n+2}^+$. 
Thus we obtain (iii).

Finally we show (iv).
If $Y$ contains $S_j^+$ and $S_j^-$ at the same time, we have $1<j<n+2$ by (iii).
Then by our choice of the components of $Y$, there exist $1\le l<l'\le m$ such that 
`$k_{j}^{(l)}<k_{j+1}^{(l)}$ and $k_{j}^{(l')}>k_{j+1}^{(l')}$' hold
(in this case $S_j^+$ and $S_j^-$ are chosen at the $l$-th and $l'$-th applications respectively) or
`$k_{j}^{(l)}>k_{j+1}^{(l)}$ and $k_{j}^{(l')}<k_{j+1}^{(l')}$' hold
(in this case $S_j^-$ and $S_j^+$ are chosen at the $l$-th and $l'$-th applications respectively),
where $(k_2^{(l)},\cdots,k_{n+2}^{(l)})$ is the sequence of $(n+1)$ integers obtained from the initial sequence $(k_1,\cdots,k_{n+2})$ by applying the procedure $(l-1)$ times.
 Suppose one of these two possibilities holds.
 Then there must exist $l''$ satisfying $l<l''<l'$ such that $k_j^{(l'')}=k_{j+1}^{(l'')}$ holds.
 But once this situation occurs, we do not select $S_j^+$ nor $S_j^-$ any more.
 Hence we always have $k_j^{(i)}=k_{j+1}^{(i)}$ for any $i\ge l''$.
 This is a contradiction.
 Thus we obtain the claim (iv).
\proofend

\begin{definition}
\label{def-m_i}
For each $1\le i\le n+2$ let $l_i$ be the multiplicity of the divisor $S_i^+$ (or $S_i^-$, equivalently) in the divisor $Y+\ol Y$, where $Y$ is the member of $|mF|$ as in Proposition \ref{prop-mt3}.
\end{definition}

Clearly we have $l_i\ge 0$, $\sum_{i=1}^{n+2}l_i=2m$ and $l_1=l_{n+2}=1$ by Prop.\ref{prop-mt3} (iii).

Let $Z, S,C_1, G_1, m$ and $Y\in|mF|$ be as in Proposition \ref{prop-mt3}.
Similarly to the notation in \cite{Hon07-4}, let $V_m$ be a subspace of $H^0(Z,mF)$ generated by the image of the natural map
\begin{equation}
H^0(F)^{G}\times H^0(F)^{G}\times\cdots\times H^0(F)^{G}\,\lra\,H^0(mF)^G
\end{equation}
given by $(s_1,s_2,\cdots,s_m)\mapsto s_1\otimes s_2\otimes\cdots\otimes s_m$.
Since $H^0(F)^{G}$ is $2$-dimensional as in Proposition \ref{prop-01}, 
$V_m$ is $(m+1)$-dimensional.
Members of the associated linear system $|V_m|$ is of the form $S_1+S_2+\cdots+S_m$, where $S_i\in |F|^G$.

\begin{prop}\label{prop-mt4}
The linear system generated by  $Y$, $\ol{Y}$, and members of the $m$-dimensional system $|V_m|$ is a $(m+2)$-dimensional subsystem (of $|mF|$).
\end{prop}

\noindent
{\bf Proof.}
Since $Y$ contains $S_1^+$ and $S_1^-$ is the only divisor that satisfies $S_1^++S_1^-\in |F|$, $Y\in|V_m|$ means $S_1^-\in |F|$.
This contradicts Proposition \ref{prop-mt3} (iii) and hence $Y\not\in |V_m|$.
Since the system $|V_m|$ is real, this implies $\ol{Y}\not\in|V_m|$.
If $\ol{Y}$ belongs to the $(m+1)$-dimensional system generated by $Y$ and $|V_m|$, then $\ol{Y}$ must contain the component $C_2$ by multiplicity $m$ since $Y|_S=mC+f-{\ol f}$, and  members of $|V_m|$ contain $C_2$ by multiplicity $m$.
On the other hand, since $\ol{Y}|_S=mC-f+{\ol f}$ and $k_2=1$, $\ol{Y}$ contains $C_2$ by multiplicity $(m-1)$.
This is a contradiction and we obtain that $Y$, $\ol{Y}$ and $|V_m|$ are linearly independent.
Hence we obtain the claim of the proposition.
\proofend

\begin{definition}\label{def-w}{\em
We denote by  $W_m$ for the linear subspace of $H^0(mF)$ generated by two sections defining the two  divisors $Y$, $\ol{Y}$, and the $(m+1)$-dimensional subspace $V_m$.
(By Proposition \ref{prop-mt4}, we have $\dim W_m=m+3$.) We denote by $\Phi_m^{G_1}:Z\to \mathbf{CP}^{m+2}$ for the meromorphic map associated to the system $|W_m|$.
}
\end{definition}

\noindent
Later we will show $W_m=H^0(mF)^{G_1}$.
We always have to keep in mind that the system $|W_m|$ is determined only after choosing a subgroup $G_1$ of $G$ (or equivalently, a component $C_1$ of $C$).

Let $\Psi_m$ be the meromorphic map associated to the $m$-dimensional system $|V_m|$ and $\Lambda_m$ a rational normal curve in $\mathbf{CP}^m$ which is the image of $\Psi_m$.
Then we have the following basic commutative diagram of meromorphic maps
\begin{equation}\label{cd1}
 \CD
Z@>{\Phi_m^{G_1}}>>\mathbf P^{\vee}W_m\\
 @V\Psi_m VV @VV{\pi_m}V\\
\Lambda_m@>>>\mathbf P^{\vee}V_m,\\
 \endCD
 \end{equation}
where $\pi_m$ is the canonical projection induced from the obvious inclusion $V_m\subset W_m$.
We note that the whole of the diagram \eqref{cd1} in not only $G_1$-equivariant but also $G$-equivariant, since the spaces $W_m$ and $V_m$ are not only $G_1$-invariant but also $G$-invariant.
As mentioned above, $W_m$ will turned out to be acted trivially by $G_1$.
On the other hand, $V_m$ is of course acted trivially by $G$.

Since the image of the restriction map $H^0(Z,mF)\supset W_m\to H^0(S,mK_S^{-1})$ is a 3-dimensional subspace generated by 3 sections defining 3 curves \eqref{gen-1}, and since the associated linear system of the latter induces the $G_1$-quotient map \eqref{quot_S} by Proposition \ref{prop-mt2}, the restriction of the meromorphic map $\Phi_m^{G_1}$ to $S$ is precisely the $G_1$-quotient map \eqref{quot_S}, whose image is a conic in a fiber of $\pi_m$ which is a  plane in $\mathbf P^{\vee}W_m$.
Thus, since $S\in|F|^G$ is an arbitrary smooth member, $\mathscr T$ can be regarded as a quotient space for the $G_1$-action on $Z$.
In view of this we introduce the following

\begin{definition}
\label{def-mt}{\em
We call the image surface $\mathscr T:=\Phi_m^{G_1}(Z)$ 
as {\em the minitwistor space of a Joyce metric with respect to the $\mathbf C^*$-subgroup $G_1$.}}
\end{definition}

In Appendix we will show that the surface $\mathscr T$ is a `canonical' quotient space of the twistor space.
Note that $\mathscr T$ is given not as an abstract complex surface but as a surface embedded in a complex projective space.
(This point will be significant later.)

Since $\Lambda_m$ is a parameter space of the pencil $|F|^G$, it follows from the diagram \eqref{cd1} that by $\pi_m$, the minitwistor space has a natural structure of a rational conic bundle over (the image of)  $\Lambda_m$.
In particular,  $\mathscr T$ is a rational surface but not a complete intersection if $m>2$ because $\Lambda_m$ is not so for $m>2$.
Moreover, the defining equation of $\mathscr T$ in $\mathbf P^{\vee}W_m$ can be explicitly determined as follows.

\begin{prop}\label{prop-mt5}
Let $Z, C_1, G_1, m, Y\in|mF|$ and $|W_m|$ be as above.
Then there exist a homogeneous coordinate $(z_0,z_1,\cdots,z_m)$ on $\mathbf P^{\vee}V_m=\mathbf{CP}^m$ and two sections $z_{m+1}, z_{m+2}\in W_m$ which satisfy the following.
(i) $\{z_0,z_1,\cdots,z_{m+2}\}$ is a basis of $W_m$.
(ii) With respect to a homogeneous coordinate $(z_0,z_1,\cdots,z_m)$ on $\mathbf{CP}^m$, the rational normal curve $\Lambda_m$ is given by
\begin{equation}\label{rnc}
\{(1,\lambda,\lambda^2,\cdots,\lambda^{m})\set\lambda\in\mathbf C\}\cup\{(0,0,\cdots,0,1)\}.
\end{equation}
(iii) $z_{m+1}$ and $z_{m+2}$ are mutually conjugate sections which define the divisors $Y$ and $\ol{Y}$ respectively.
(iv) The minitwistor space $\mathscr T=\Phi_m^{G_1}(Z)$ satisfies not only equations in the defining ideal of $\Lambda_m$ in $\mathbf{CP}^m$ but also a quadratic equation of the following form
\begin{equation}\label{eqn_mt}
z_{m+1}z_{m+2}=Q(z_0,z_1,\cdots,z_m),
\end{equation}
where $Q$ is a quadratic homogeneous polynomial with real coefficients.
(v) The degree of the surface $\mathscr T$ in $\mathbf{CP}^{m+2}$ is $2m$.
\end{prop}

\noindent {\bf Proof.}
This can be proved by an argument  analogous to that of  \cite[Theorem 2.11]{Hon07-4}. So here we only sketch a proof.
For $1\le i\le n+2$ let $e_i$ be a section of a line bundle $[S_i^+]$ which defines $S_i^+$, and put $u_i=e_i\otimes \ol{e}_i\in H^0(F)$.
We adopt $\{u_1,u_{n+2}\}$ as a basis of $H^0(F)^G$.
Let $0=\lambda_1,\lambda_2,\cdots,\lambda_{m+1}$ be a sequence of real numbers satisfying $u_i=u_1-\lambda_iu_{n+2}$ and put $\lambda_{m+2}=\infty(=-\infty)$.
Then either $\lambda_1<\lambda_2<\cdots<\lambda_{m+2}=\infty$ or 
$\lambda_1>\lambda_2>\cdots>\lambda_{m+2}=-\infty$ holds.
For $0\le i\le m$ we put $z_i=u_1^iu_{n+2}^{m-i}$ and use $\{z_0,z_1,\cdots,z_m\}$ as a basis of $V_m$.
Then the rational normal curve $\Lambda_m$ is written in the form \eqref{rnc}.

Since the divisor $Y\in |mF|$ is a sum of degree one divisors as in Proposition \ref{prop-mt3}, its defining section is a product of $e_i$ and $\ol{e}_i$ $(1\le i\le n+2)$, where $e_i$ or $\ol{e}_i$ is selected iff $S_i^+$ or $S_i^-$ is contained as an irreducible component of $Y$ respectively.
(If $Y$ contains $S_i^+$ by multiplicity $l_i$, then $e_i$ is selected $l_i$ times.)
Let $z_{m+1}$ be a defining section of $Y$ obtained in this way, and put $z_{m+2}=\ol{z}_{m+1}$.
Then by Proposition \ref{prop-mt4} and Definition \ref{def-w}, $\{z_0,z_1,\cdots,z_{m+2}\}$ is a basis of $W_m$.
Moreover, by Definition \ref{def-m_i}, we have
\begin{align}
z_{m+1}z_{m+2}&=c\, (e_1\ol{e}_1)^{l_1} (e_2\ol{e}_2)^{l_2}\cdots(e_{n+2}\ol{e}_{n+2})^{l_{n+2}}\label{key1}\\
&=c\,u_1u_{2}^{l_2}u_3^{l_3}\cdots u_{n+1}^{l_{n+1}}u_{n+2}\\
&=c\,u_1(u_1-\lambda_2u_{n+2})^{l_2}(u_1-\lambda_3u_{n+2})^{l_3}\cdots (u_1-\lambda_{n+1}u_{n+2})^{l_{n+1}}u_{n+2},
\label{dcp1}
\end{align} 
where $c$ is a non-zero real constant which can be taken as either $1$ or $-1$.
Then expanding \eqref{dcp1} and using the relation $l_1+l_2+\cdots+l_{n+2}=2m$, we obtain that the right-hand side of \eqref{dcp1} can be written as a quadratic polynomial $Q$ of $z_0,z_1,\cdots,z_m$ with real coefficients. 
(Of course, $Q$ is determined only up to the defining ideal of $\Lambda_m$. 
But it is explicitly computable since $l_1,\cdots,l_{n+2}$ are computable.)
Thus we obtain (iv).

Finally (v) is immediate since by (iv) $\mathscr T$ is an intersection of the rational normal scroll of 2-planes parameterized by $\Lambda_m$ whose degree is $m$, and the quadratic hypersurface \eqref{eqn_mt}.
\proofend

\bigskip
We note that $\{\lambda_1,\cdots,\lambda_{n+2}\}$ in the proof is nothing but the conformal invariant  of Joyce metrics.
Note also that $l_2,\cdots,l_{n+2}$ in \eqref{dcp1} are readily computable numbers.

The two subspaces $V_m,W_m\subset H^0(Z,mF)$ have the following intrinsic characterizations.

\begin{prop}\label{prop-char}
(i) $V_m=H^0(Z,mF)^G$ holds.
(ii) $W_m=H^0(Z,mF)^{G_1}$ holds.	
\end{prop}

\noindent
{\bf Proof.}
For (i), the inclusion $V_m\subset H^0(Z,mF)^G$ is obvious.
Suppose $\dim H^0(Z,mF)^G=(m+1)+\alpha$, $\alpha>0$.
Then the image of the meromorphic map associated to $|mF|^G$ must still be 1-dimensional, since general orbits of $G$-action on $Z$ are 2-dimensional.
Let $\Lambda_m'$ be the image, which is a non-degenerate curve in $\mathbf{CP}^{m+\alpha}$. 
Considering the natural projection $\mathbf{CP}^{m+\alpha}\to\mathbf{CP}^m$ induced by the inclusion $V_m\subset H^0(Z,mF)^G$, the map $\Psi_m$ is factorized as $Z\to\Lambda'_m\to\Lambda_m$.
Here, since $\alpha>0$ and $\Lambda_m'$ is non-degenerate, the degree of the map $\Lambda'_m\to\Lambda_m$ is at least two.
On the other hand, general fiber of $\Psi_m:Z\to\Lambda_m$ is connected.
This is a contradiction and hence we obtain $\alpha=0$.
Thus we obtain (i).

For (ii), identifying $\Lambda_m$ and its image into $\mathbf P^{\vee}V_m$, by
 \eqref{cd1}, we obtain the following commutative diagram of meromorphic maps
\begin{equation}\label{cd2}
\xymatrix{
   Z \ar@{->}[r]^{{\Phi_m^{G_1}}}\ar@{->}[d]_{\Psi_m}  & \mathscr T  \ar@{->}[dl]^{\pi_m}\\
   \Lambda_m. & \\
}
 \end{equation}
If the defining section $z_{m+1}$ of $Y$ is not $G_1$-invariant, $G_1$ acts non-trivially on the image surface $\mathscr T$ (since $\mathscr T$ is non-degenerate in $\mathbf{CP}^{m+2}$).
Since $G$ acts trivially on $\Lambda_m$, by \eqref{cd2}, this means that $G_1$ acts non-trivially on fibers of $\mathscr T\to\Lambda_m$.
This is a contradiction since the restriction $\Phi_m^{G_1}|_S$ can be identified with the $G_1$-quotient map \eqref{quot_S}.
Therefore the section $z_{m+1}$ is also $G_1$-invariant.
Then by reality, $z_{m+2}$ is also $G_1$-invariant. 
Hence, since $V_m=H^0(mF)^G\subset H^0(mF)^{G_1}$,
we obtain $W_m\subset H^0(Z, mF)^{G_1}$. 
Suppose that $\dim H^0(mF)^{G_1}=(m+3)+\alpha$, $\alpha>0$.
Then by the same argument for $\Psi_m$ given above, the map $\Phi_m^{G_1}$ is factorized as $Z\to\mathscr T'\to\mathscr T$, where $\mathscr T'$ is a non-degenerate surface in $\mathbf{CP}^{m+2+\alpha}$, and the map $\mathscr T'\to\mathscr T$ is generically $l:1$ with $l\ge 2$.
This contradicts the fact that general fibers of $\Phi_m^{G_1}$ are 
the closures of  $G_1$-orbits that are necessarily irreducible.
Thus we obtain $\alpha=0$, as required.
\proofend


\bigskip
The proposition means that the rational map $\pi_m$ can also be regarded as a quotient map for $G/G_1\,(\simeq\mathbf C^*)$-action on $\mathscr T$.
As to the reducible fibers of $\pi_m$ (which are necessarily a sum of two lines), the equations \eqref{key1}--\eqref{dcp1} directly imply the following, which means that the complex structure of $\mathscr T$ depends on the conformal invariants of Joyce metrics.

\begin{prop}\label{prop-rf1}
The rational map $\pi_m:\mathscr T\to\Lambda_m$ has reducible fibers precisely over the points $(1,\lambda_i,\lambda_i^2,\cdots,\lambda_i^m)\in\Lambda_m$ and $(0,0,\cdots,0,1)\in\Lambda_m$ in the coordinate of Proposition \ref{prop-mt5}, where the index $i$ of $\lambda_i$ must satisfy $l_i>0$.
\end{prop}

Note that if  $m\ge 2$ (i.e.\,if $G_1$ does not act semi-freely on $Z$), the complex structure of $\mathscr T$ actually varies in accord with the conformal invariant of Joyce metrics
since the number of reducible fibers of $\pi_m$ is greater than $3$ by Proposition \ref{prop-rf1}.
Thus our minitwistor spaces constitute a non-trivial moduli space, if $m\ge 2$.
The dimension of the moduli space is given by
\begin{equation}\label{rf2}
\#\left\{1\le i\le n+2\set l_i>0
\right\}-3,
\end{equation}
where the first term is the number of reducible fibers of $\pi_m$ and $3$ is the dimension of $PGL(2,\mathbf C)$ acting on $\Lambda_m\simeq\mathbf{CP}^1$.
Proposition \ref{prop-rf1} also means that, if $l_i=0$, the fiber over the point $(1,\lambda_i,\lambda_i^2,\cdots,\lambda_i^m)$ is irreducible. (The condition $l_i=0$ is equivalent to `$S_i^+\not\subset Y$ and $S_i^-\not\subset Y$'.)
For these $i$, the twistor line $L_i=S_i^+\cap S_i^-$ has the following important property.

\begin{prop}\label{prop-irf}
If $l_i=0$, then  $L_i \subset Z^{G_1}$ holds.
\end{prop}

\noindent
{\bf Proof.}
By Proposition \ref{prop-rf1} the fiber $\pi_m^{-1}(1,\lambda_i,\lambda_i^2,\cdots,\lambda_i^m)$ is irreducible  if $l_i=0$.
Since Bs\,$|W_m|\subset C$ and $L_i\not\subset C$, the image $\Phi_m^{G_1}(L_i)$ makes sense for any $1\le i\le n+2$ and they are $G$-invariant.
By the diagram \eqref{cd2},  $\Phi_m^{G_1}(L_i)\subset \pi_m^{-1}(1,\lambda_i,\lambda_i^2,\cdots,\lambda_i^m)$.
If $\Phi_m^{G_1}(L_i)$ is a point, the point becomes $G$-fixed real point of $\mathscr T$. 
But since the fiber $\pi_m^{-1}(1,\lambda_i,\lambda_i^2,\cdots,\lambda_i^m)$ is irreducible, there exists no such a point.
Hence  $\Phi_m^{G_1}(L_i)$ is a curve, coinciding with the fiber itself.
Since  $\Phi_m^{G_1}$ is $G$-equivariant and $G_1$ acts trivially on $\mathscr T$, $L_i$ is acted trivially by $G_1$, as desired.
\proofend

\bigskip


%

%

Next we determine the singularities of the  minitwistor spaces $\mathscr T$.
For this we first note that the center (= the indeterminacy locus) of the projection $\pi_m:\mathbf P^{\vee}W_m=\mathbf{CP}^{m+2}\to\mathbf P^{\vee} V_m=\mathbf{CP}^m$ in the diagram \eqref{cd1} is a line explicitly given by $\{z_0=z_1=\cdots=z_m=0\}$ in the coordinate of Proposition \ref{prop-mt5}.
Let $l_{\infty}$ be this line.
Since fibers of $\pi_m$ are exactly  2-planes which contain the line $l_{\infty}$, the rational scroll $\pi_m^{-1}(\Lambda_m)$ contains $l_{\infty}$.
Further,   $\pi_m^{-1}(\Lambda_m)$ has cyclic quotient singularities along $l_{\infty}$, where the singularity looks like $\mathbf C^2/\mathbf Z_m$, with the generator of  $\mathbf Z_m$ acting on $\mathbf C^2$ as $(z,w)\mapsto (\zeta z,\zeta w)$ with $\zeta=e^{\frac{2\pi i}{m}}$.
(So it is not a rational double point if $m>2$.)
On the other hand, the quadratic hypersurface \eqref{eqn_mt} intersects $l_{\infty}$ at 2 points $P_{\infty}:=(0,\cdots,0,1,0)$ and $\ol{P}_{\infty}:=(0,\cdots,0,0,1)$.
This means that the indeterminacy locus of the projection $\pi_m:\mathscr T\to\Lambda_m$ consists of  $\{P_{\infty},\ol{P}_{\infty}\}$, and these are also  singular points of the surface $\mathscr T$.
Further, since the intersection number of the quadratic \eqref{eqn_mt} and $l_{\infty}$ is of course two, the intersection must be transversal.
Hence $P_{\infty}$ and $\ol{P}_{\infty}$ are cyclic quotient singularities of $\mathscr T$ of the above kind.

If we blow-up $\mathbf{CP}^{m+2}$ along $l_{\infty}$, the indeterminacy locus is eliminated and $\mathbf{CP}^{m+2}$ is transformed into the total space of the $\mathbf{CP}^2$-bundle $\mathbf P(\mathscr O(1)^{\oplus 2}\oplus\mathscr O)$ over $\mathbf{CP}^m$.
Let $\hat{\mathscr T}$ be the strict transform of $\mathscr T$ under this blowing-up.
Then the natural map $\hat{\mathscr T}\to\mathscr T$ obviously resolves the two   singularities $P_{\infty}$ and $\ol{P}_{\infty}$.
We denote by $\Gamma$ and $\ol{\Gamma}$ for the exceptional curves over $P_{\infty}$ and $\ol{P}_{\infty}$ respectively.

Since $\mathscr T$ is contained in $\pi^{-1}_m(\Lambda_m)$, $\hat{\mathscr T}$ is contained in the restriction of the $\mathbf{CP}^2$-bundle
$\mathbf P(\mathscr O(1)^{\oplus 2}\oplus\mathscr O)\to\mathbf{CP}^m$ onto the curve $\Lambda_m$.
Since the degree of the curve $\Lambda_m$ in $\mathbf{CP}^m$ is $m$, this restriction is identified with the $\mathbf{CP}^2$-bundle 
$\mathbf P(\mathscr O(m)^{\oplus 2}\oplus\mathscr O)\to\Lambda_m\simeq\mathbf{CP}^1$.
Thus the surface $\hat{\mathscr T}$ is  embedded in this $\mathbf{CP}^2$-bundle.
By the equations \eqref{key1}--\eqref{dcp1}, the defining equation of $\hat{\mathscr T}$ in this bundle is explicitly given by
\begin{equation}
\xi_{1}\xi_{2}
=c\,u(u-\lambda_2)^{l_2}(u-\lambda_3)^{l_3}\cdots(u-\lambda_{n+1})^{l_{n+1}}.
\label{dcp2}
\end{equation} 
where $u=u_1/u_{n+2}$ is a non-homogeneous coordinate on $\Lambda_m$ and $(\xi_1,\xi_2)=(z_{m+1}/u_{n+2}^{m},z_{m+2}/u_{n+2}^m)$ represents points on the vector bundle $\mathscr O(m)^{\oplus 2}\to\Lambda_m$.
From this expression, we can readily deduce the following.

\begin{prop}\label{prop-sing_T}
Let $\mathscr T$ be the minitwistor space of the twistor space of a Joyce metric on $n\mathbf{CP}^2$ with respect to a $\mathbf C^*$-subgroup $G_1$ of $G$ as above.
Then the singular locus of $\mathscr T$ consists of the following.
(a) The conjugate pair $P_{\infty}$ and $\ol{P}_{\infty}$ which are isomorphic to a cyclic quotient singularity $\mathbf C^2/\mathbf Z_m$ as above,
(b) The real points of the form $(1,\lambda_i,\lambda_i^2,\cdots,\lambda_i^m,0,0)\in\mathscr T\subset\mathbf{CP}^{m+2}$ where 
the index $i$ satisfies $l_i>1$.
Further, this singular point is $A_{l_i-1}$-singularity of $\mathscr T$.
\end{prop}

\section{Projective models of the twistor spaces of Joyce metrics}
In the previous section, we have explicitly given a linear subsystem $|W_m|$ of $|mF|$ which induces a $G_1$-quotient map $\Phi_m^{G_1}:Z\to\mathscr T$ for the twistor spaces of Joyce metrics.
The quotient surface $\mathscr T$, which we call the minitwistor space with respect to $G_1$, was shown to be a rational surface with  mild singularities.
In this section, by investigating the meromorphic map $\Phi_m^{G_1}$ more in detail, we realize a projective model of the twistor space as a conic bundle over (a resolution of) the minitwistor space.

For this purpose, we recall that
$|W_m|$ is $(m+2)$-dimensional system generated by the $m$-dimensional subsystem $|V_m|$ (composed of the pencil $|F|^G$) and two divisors $Y$ and $\ol Y$ which are explicitly given by \eqref{Y} and its conjugation.
We rewrite \eqref{Y} as
\begin{align}\label{Y3}
Y=\sum_{1\le i\le n+2}(l_i^+S_i^++l_i^-S_i^-).
\end{align}
Then by Proposition \ref{prop-mt3} (iv), at least one of $l_i^+=0$ and $l_i^-=0$ holds for each $i$.
Further, by Proposition \ref{prop-mt3} (iii), $l_1^+=l_{n+2}^-=1$ and $l_1^-=l_{n+2}^+=0$.
We have $l_i=l_i^++l_i^-$ (see Definition \ref{def-m_i}). 
We also have
\begin{align}\label{Y4}
m=\sum_{1\le i\le n+2}l_i^+=\sum_{1\le i\le n+2}l_i^-
\end{align}
which immediately follows from the expression \eqref{Y}.
In the following argument we frequently use these numbers.

\subsection{Partial elimination of  the base locus of the system  $|W_m|$ }
Because the base locus of a linear system is the intersection of all divisors generating  the system, and since in the present situation the base locus of the subsystem $|V_m|$ of $|W_m|$ is exactly the cycle $C$, we have Bs\,$|W_m|=Y\cap \ol{Y}\cap C$.
Further, since $Y$ contains the component $S_{n+2}^-$ containing $C_1$ and also $S_1^+$ containing $\ol{C}_1$, we have $Y\supset C_1\cup\ol C_1$.
Hence we also have $\ol Y\supset C_1\cup\ol C_1$.
These mean Bs\,$|W_m|\supset C_1\cup\ol C_1$.
Since the closure of general $G_1$-orbits go through $C_1$ and $\ol C_1$, $C_1\cup\ol C_1$ is, in some sense, the `principal part' of the base locus.
In the following, we give a sequence of blowing-ups which eliminates the base locus  that intersects $C_1\cup\ol{C}_1$.

Before doing so, we give some remark.
Since we know generating divisors of $|W_m|$ and the way how they intersect, it is in principle possible to give a full elimination of the base locus, which may lead us to an explicit construction of the twistor spaces of arbitrary Joyce metrics.
Actually, for the twistor space of LeBrun metric with $K$-action, if we take the subgroup $G_1$  acting semi-freely on $n\mathbf{CP}^2$, we have $m=1$ and Bs\,$|W_1|$  is eliminated by just blowing-up  $C_1\cup\ol{C}_1$.
This seems to be a basic reason why the twistor spaces of LeBrun metrics can be constructed relatively easily.
But for other twistor spaces of Joyce metrics, the required blow-ups are so complicated and it looks not easy to accomplish. 
So here we only  give a partial elimination, which is enough for the purpose of giving a projective model.

As the first step 
let $\mu_1:Z_1\to Z$ be the blowing-up along $C_1\cup\ol{C}_1$, and $E_1$ and $\ol{E}_1$ the exceptional divisors over $C_1$ and $\ol{C}_1$ respectively.
The two divisors $S_1^-$ and $S_{n+2}^-$ contain $C_1$ and intersect transversally along $C_1$.
Further, if we put $l=(C_1)_S^2$ for the self-intersection number in $S$ ($S$ is a smooth member of $ |F|^G$ as before), we have $(C_1)^2_{S_1^-}=(C_1)^2_{S_{n+2}^-}=l+1$ (cf.\,\cite[p.\,241 (13)]{F00}).
Therefore $N_{C_1/Z}$  is isomorphic to $\mathscr O(l+1)^{\oplus 2}$.
Hence $E_1$ (and $\ol E_1$) is isomorphic to $\mathbf{CP}^1\times\mathbf{CP}^1$.
Further, on $Z_1$, the intersections $S_1^-\cap E_1$ and $S_{n+2}^-\cap E_1$ are disjoint sections of the natural projection $E_1\to C_1$ whose self-intersection numbers (in $E_1$) are zero.
On the other hand, since both $S_1^+$ and $S_{n+2}^+$ intersect $C_1$ transversally at a unique point respectively, the intersections $S_1^+\cap E_1$ and $S_{n+2}^+\cap E_1$ are (different) fibers of  $E_1\to C_1$.
Thus the restrictions $S_1^{\pm}\cap E_1$ and $S_{n+2}^{\pm}\cap E_1$ form a `square' in $E_1\simeq\mathbf{CP}^1\times\mathbf{CP}^1$.

Since Bs\,$|F|^G=C\supset C_1\cup\ol C_1$, the strict transform of the pencil $|F|^G$ makes sense, and we still denote it by $|F|^G$.
This  determines a pencil on $E_1$ and $\ol E_1$ by restriction.
Clearly we can take $(S_1^++S_1^-)\cap E_1$ and $(S_{n+2}^++ S_{n+2}^-)\cap E_1$ as generators of the  pencil on $E_1$.
It follows that  this pencil on $E_1$ has bidegree $(1,1)$ and has precisely two base points.
As long as $2\le i\le n+1$, the intersections $S_i^-\cap E_1$ are irreducible members of this pencil.
On the other hand, since $S_i^+\cap C_1=\emptyset$ for $2\le i\le n+1$ 
we have $S_i^+\cap E_1=\emptyset$ for these $i$.
Since the pencil $|F|^G$ is real, the analogous result holds for the  restriction on $\ol E_1$.

With these information, we investigate transformations of the system $|W_m|$ on $Z$ into $Z_1$.
Members of the $m$-dimensional subsystem $|V_m|$ clearly contain the curves $C_1$ and $\ol{C}_1$ by multiplicity exactly $m$.
Hence the set $\{\mu_1^{-1}(D)-m(E_1+\ol{E}_1)\set D\in |V_m|\}$ of divisors on $Z_1$ is still a linear system,
 still denoted by $|V_m|$.
Further, by \eqref{Y4},  both $Y$ and $\ol{Y}$ also contain $C_1$ and $\ol{C}_1$ by multiplicity exactly $m$.
Hence the set $\{\mu_1^{-1}(D)-m(E_1+\ol{E}_1)\set D\in |W_m|\}$  is also a linear system,
 still denoted by $|W_m|$.
Then  $|W_m|$ (on $Z_1$) contains  $|V_m|$ as a  subsystem and is still generated by $|V_m|$, $Y$ and $\ol Y$, where $Y$ remains expressed as \eqref{Y3} if we interpret  all the divisors in the right-hand side as being those on $Z_1$.

As in the case on $Z$, we have Bs\,$|W_m|=$\,Bs\,$|V_m|\cap Y\cap \ol Y$ on $Z_1$.
We first determine Bs\,$|V_m|$ (on $Z_1$).
Since  $m(S_1^++S_1^-)\in |V_m|$ and $m(S_{n+2}^++S_{n+2}^-)\in |V_m|$ still hold,
we have Bs\,$|V_m|\subset (S_1^+\cup S_1^-)\cap(S_{n+2}^+\cup S_{n+2}^-)$.
The latter is $(S_1^+\cap S_{n+2}^+)\cup(S_1^+\cap S_{n+2}^-)\cup
( S_1^-\cap S_{n+2}^+)\cup( S_1^-\cap S_{n+2}^-)$.
Since $S_1^+\cap S_{n+2}^+=\ol C_1$ on $Z$ and $(S_1^+\cap \ol E_1)\cap (S_{n+2}^+\cap\ol E_1)=\emptyset$ on $Z_1$ as noted above, $S_1^+\cap S_{n+2}^+=\emptyset$ on $Z_1$.
Hence $ S_1^-\cap S_{n+2}^-=\emptyset $ on $Z_1$.
Further since $S_1^+\cap S_{n+2}^-=C_2\cup C_3\cup\cdots \cup C_{n+2}$ on $Z$, the same equality holds also on $Z_1$.
Hence we have $S_1^-\cap S_{n+2}^+=\ol C_2\cup \ol C_3\cup\cdots \cup\ol C_{n+2}$ on $Z_1$.
From these we obtain, on $Z_1$,
\begin{align}\label{bs5}
{\rm{Bs}}\,|V_m|=(C_2\cup C_3\cup\cdots \cup C_{n+2})\cup(\ol C_2\cup \ol C_3\cup\cdots \cup\ol C_{n+2}).
\end{align}

Now suppose that $m=1$.
(This is equivalent, by Proposition \ref{prop-char-LB}, to the assumption that $Z$ is the twistor space of LeBrun metric with $K$-action, and the $G_1$ is acting semi-freely on $Z$.)
Then since $k_i=1$ for all $2\le i\le n+1$,
we have $Y=S_1^++S_{n+2}^-$ and $\ol Y=S_1^-+S_{n+2}^+$ by \eqref{Y}.
Since $(S_1^+\cup S_{n+2}^-)\cap (S_1^-\cup S_{n+2}^+)=(S_1^+\cap S_1^-)\cup(S_1^+\cap S_{n+2}^+)\cup(S_{n+2}^-\cap S_1^-)\cup(S_{n+2}^-\cap S_{n+2}^+)$, $S_1^+\cap S_1^-=L_1$,
$S_{n+2}^-\cap S_{n+2}^+=L_{n+2}$, and $S_1^+\cap S_{n+2}^+=S_{n+2}^-\cap S_1^-=\emptyset$ on $Z_1$, we obtain $Y\cap \ol Y\subset L_1\cup L_{n+2}$.
But $L_1$ and $ L_{n+2}$ are  disjoint from the right-hand side of \eqref{bs5}, at least on $Z_1$.
Hence we obtain Bs\,$|W_m|=\emptyset$.
Namely if $m=1$, Bs\,$|W_m|$ on $Z$ is completely eliminated if we just blow-up $C_1\cup\ol C_1$.

In the following we suppose $m\ge 2$.
By \eqref{bs5}, members of $|V_m|$ (on $Z_1$) actually contain the 4 curves $C_2,\ol C_{n+2},\ol C_2$ and $ C_{n+2}$, and the  multiplicities are $m$.
Further, by Proposition \ref{prop-mt3} (ii) the member $Y$ contains  $C_2$ and $\ol{C}_{n+2}$ by multiplicity $(m+1)$ and $(m-1)$ respectively, and $\ol{Y}$ contains  $C_2$ and $\ol{C}_{n+2}$ by multiplicity $(m-1)$ and $(m+1)$ respectively.
These mean  Bs\,$|W_m|\supset C_2\cup\ol C_{n+2}\cup\ol C_2\cup C_{n+2}$.


So
let $\mu_2:Z_2\to Z_1$ be the blowing-up along  $C_2\cup\ol{C}_{n+2}\cup\ol{C}_2\cup C_{n+2}$, and $E_2,\ol{E}_{n+2}, \ol{E}_2$ and $E_{n+2}$ the exceptional divisors respectively.
By the above investigation of Bs\,$|W_m|\subset Z_1$, the pull-back $\mu_2^*|W_m|$  has $(m-1)(E_2+\ol{E}_{n+2}+ \ol{E}_2+E_{n+2})$ as its fixed components.
We still denote by $|W_m|$ for the $(m+2)$-dimensional linear system on $Z_2$ whose members are the total transforms of the system $|W_m|$ on $Z_1$ with these fixed components subtracted.
Also, we still denote by $|V_m|$ for the $m$-dimensional linear system on $Z_2$ whose members are the total transforms of the system $|V_m|$ on $Z_1$ with the same fixed components   subtracted.

By the same argument for the system $|V_m|$ on $Z_1$, we have, on $Z_2$, 
the linear system $|V_m|$  is generated by the following $(m+1)$ divisors:
\begin{equation}\label{rest-6}
k(S_1^++S_1^-)+(m-k)(S_{n+2}^++S_{n+2}^-)+(E_2+\ol{E}_{n+2}+\ol E_2+E_{n+2}),\,\,\,\,0\le k\le m,
\end{equation}
where the last term comes from the fact that all members of $|V_m|$ in $Z_1$ contain $C_2,\ol{C}_{n+2}, \ol C_2$ and $C_{n+2}$ by multiplicity precisely $m$ which exceeds $(m-1)$ by one.
By considering the two divisors obtained by putting $k=0$ and $k=m$ and taking their intersection, we obtain 
\begin{align}\label{rest-22}
{\rm{Bs}}\,|V_m|=\{(C_3\cup C_4\cup\cdots \cup C_{n+1})\cup(\ol C_3\cup \ol C_4\cup\cdots \cup\ol C_{n+1})\}\cup E_2\cup\ol E_2\cup E_{n+2}\cup\ol E_{n+2}.
\end{align}
Of course, the 4 exceptional divisors in \eqref{rest-22} are disjoint.
Since we are eliminating the base locus of $|W_m|$ which intersects $E_1\cup\ol E_1$, we ignore the curve part of \eqref{rest-22}.
For another generator $Y$, we have
\begin{equation}\label{rest-4}
Y=S_1^++\sum_{2\le i\le n+1}(l_i^+S_i^++l_i^-S_{i}^-)+S_{n+2}^-+2(E_2+E_{n+2}),
\end{equation}
where the last term is a consequence of the fact that $Y$ contains the curves $C_2$ and $C_{n+2}$ by multiplicity $(m+1)$ which exceeds $(m-1)$ by two.
By conjugation, we have 
\begin{equation}\label{rest-5}
\ol Y=S_1^-+\sum_{2\le i\le n+1}(l_i^+S_i^-+l_i^-S_{i}^+)+S_{n+2}^++2(\ol E_2+\ol E_{n+2}).
\end{equation}
We are determining the intersection $(E_2\cup\ol E_2\cup E_{n+2}\cup\ol E_{n+2})\cap Y\cap\ol Y$.
Since $E_2\subset Y$, we have $E_2\cap Y\cap\ol Y=E_2\cap\ol Y$.
Further, we have $E_2\cap\ol E_2=E_2\cap \ol E_{n+2}=\emptyset$, and also $E_2\cap S_i^+=\emptyset$ if $3\le i\le n+2$ (since $C_2\cap S_i^+=\emptyset$ on $Z$ for these $i$).
On the other hand, $E_2\cap S_2^+$ is a curve since $C_2$ intersects $S_2^+$ transversally at a point on $Z$ (and on $Z_1$).
But since this curve is clearly disjoint from $E_1$, we ignore.
Furthermore, $E_2\cap S_i^-$ is a curve intersecting $E_1$ transversally at a point, for $2\le i\le n+2$.
Also, $E_2\cap S_1^-=\emptyset$ since  $C_2$ and $S_1^-$ are made apart by the first blow-up $\mu_1$.
Thus we obtain that, among components of $E_2\cap Y\cap\ol Y$, only the curves $S_i^-\cap E_2$, $2\le i\le n+1$ intersect $F_1$ (for $i$ satisfying $l_i^+>0$, of course).
By a similar consideration, we deduce that the irreducible components of $\ol E_{n+2}\cap Y\cap\ol Y$ which intersects $E_1$ are precisely given by
$S_i^-\cap \ol E_{n+2}$ where $i$ satisfies $2\le i\le n+1$ and $l_i^->0$.
Thus we have obtained that, on $Z_2$, the components of Bs\,$|W_m|$ which intersect $E_1$ are explicitly given by
\begin{equation}\label{bc1}
\left\{E_2\cap S_i^-\set 1<i<n+2,\,l_i^+>0\right\}\,\bigcup\,
\left\{\ol{E}_{n+2}\cap S_i^-
\set 1<i<n+2,\,l_i^->0\right\}
 \end{equation}
All these are clearly smooth rational curves, intersecting $E_1$ transversally at a point respectively.
Members of the system $|V_m|$ on $Z_2$ contain all these curves by multiplicity one,
since the coefficients of $E_2$ and $\ol{E}_{n+2}$ in \eqref{rest-6} are one.
On the other hand, the divisor $Y$ on $Z_2$ contains the curve $E_2\cap S_i^-$ by multiplicity two 
(because of the coefficient of $E_2$ in \eqref{rest-4})
and the curve $\ol{E}_{n+2}\cap S_i^-$ by multiplicity $l_i^-$
(because of the coefficient of $S_i^-$ in \eqref{rest-4}).
Similarly, the divisor $\ol{Y}$ on $Z_2$ contains the curve $E_2\cap S_i^-$ by multiplicity $l_i^+$ and the curve $\ol{E}_{n+2}\cap S_i^-$ by multiplicity two.
Thus we have  obtained the base locus of the system $|W_m|$ on $Z_2$ intersecting the divisor $E_1$, and also  multiplicities of the generators  along the base curves.
By considering the conjugations, we obtain  the base locus of the system $|W_m|$ on $Z_2$ which intersects the conjugate divisor $\ol{E}_1$ as well as the multiplicities of generators, in a concrete form.

As the third step for the elimination of the base locus, let $\mu_3:Z_3\to Z_2$ be the blowing-up along the curves \eqref{bc1} and their conjugations, and $F_i$ the exceptional divisor over the curve 
 $E_2\cap S_i^-$ or $\ol{E}_{n+2}\cap S_i^-$ appearing in \eqref{bc1}.
 (Note that for each $i$, at most one of $l_i^+>0$ and $l_i^->0$ holds.
 Hence this definition makes sense.)
 Then by the above consideration of the multiplicities,  the fixed components of $\mu_3^*|W_m|$ is the sum of all exceptional divisors $F_i$ and $\ol{F}_i$ of $\mu_3$, where their multiplicities are exactly one for every these divisors.
 So as the transformation of the system $|W_m|$  into $Z_3$, we  consider the system whose members are the total transforms of the members of $|W_m|$ with the divisor $\sum (F_i+\ol{F}_i)$ removed, where $i$ runs  satisfying $l_i>0$.
 We again denote $|W_m|$ for this linear system on $Z_3$.

 By  computations similar to those on $Z_1$ and $Z_2$  above, we deduce that, on $Z_3$, the components of Bs\,$|W_m|$ which intersect $E_1$ are given by
 \begin{align}\label{bs55}
 \left( \bigcup_i F_i\cap E_2
\right)
\cup
\left( \bigcup_i F_i\cap \ol{E}_{n+2}
\right)
 \end{align}
 where $i$ runs  satisfying $l_i^+\ge 2$ in the first union and $l_i^-\ge 2$ in the second union.
 Further, the divisor $Y$ (on $Z_3$) contain the curves $F_i\cap \ol{E}_{n+2}$ and $F_i\cap E_2$ by multiplicities $(l_i^--1)$ and $1$ respectively.
The member $\ol Y$ contains these curves by multiplicities $1$ and $(l_i^+-1)$ respectively.
On the other hand, 
 members of $|V_m|$ (on $Z_3$) contains these curves by multiplicity one.

In particular, if $l_i=1$ for any $1\le i\le n+2$ (which is actually the case for some Joyce metrics and some appropriate subgroup $G_1$; see \S 5.3),  \eqref{bs55} is empty set and we conclude that the linear system $|W_m|$ on $Z_3$ is base point free, at least in a neighborhood of $E_1\cup\ol{E}_1$.
 If $l_i>1$ for some $1<i<n+2$, 
we take further blowing-up $\mu_4:Z_4\to Z_3$ along the curves \eqref{bs55}
and their conjugations.
Let  $F_i'$ and $\ol{F}'_i$ be the exceptional divisors over `$F_i\cap \ol{E}_{n+2}$ or $F_i\cap E_2$' and `$\ol F_i\cap \ol{E}_{n+2}$ or $\ol F_i\cap E_2$' respectively.
(As before, precisely one of $F_i\cap \ol{E}_{n+2}\neq\emptyset$ or
$F_i\cap E_2\neq\emptyset$ occurs.)
 We take $\mu_4^*|W_m|$  and subtract its fixed components $\sum(F_i'+\ol{F}'_i)$, where $i$ runs  satisfying $l_i>2$.
Then the base locus of the resulting linear  system on $Z_4$ intersecting $E_1$ is given by the union of rational curves $F_i'\cap \ol{E}_{n+2}$ and $F'_i\cap E_2$ where $i$ satisfies $l_i>2$ this time.
In particular, if $l_i\le 2$ for all $i$, the linear system on $Z_4$ is base point free, at least in a neighborhood of $E_1\cup\ol{E}_1$.
If $l_i>2$ for some $i$, we need further blowing-up along $F'_i\cap E_2$ and $F'_i\cap \ol{E}_{n+2}$  where $i$ satisfies $l_i>2$.

Since the multiplicities along the base curves decrease by one  for each blowing-ups as before,
this blowing-up process stops in a finite number of times.
(Namely, if $M:=\max\{l_i\set 1\le i\le n+2\}$, the linear system on $Z_{M+2}$ becomes base point free in a neighborhood of $E_1\cup\ol{E}_1$.)
Consequently, the base locus of the linear system $|W_m|$ on $Z$ is removed by this process,  at least in a neighborhood of $C_1\cup\ol{C}_1$.
Note that these blowing-ups are explicitly and uniquely specified.
Let $\tilde{Z}$ (= $Z_{M+2}$ in the above notation) be the resulting manifold.
Since all the centers of the blow-ups are $G$-invariant, $\tilde{Z}$ has a natural $G$-action.
We denote by $|\tilde{W}_m|$ for the linear system on $\tilde{Z}$ corresponding to the system $|W_m|$ on $Z$.
Then the meromorphic map associated to the system $|\tilde{W}_m|$ is precisely  the composition of the sequence of  blowing-ups  $\tilde{Z}\to\cdots\to Z_1\to Z$ and $\Phi_m^{G_1}:Z\to \mathscr T$.
 Of course, this meromorphic map is a morphism in a neighborhood of $E_1\cup\ol{E}_1$.
We note that our blow-ups also eliminate all the base curves that are contained in the 4 divisors $E_2,\ol{E}_{n+2},\ol{E}_2$ and $E_{n+2}$.
Therefore, identifying the curves $C_i$ and $\ol{C}_i$ in $Z$ with their strict transforms into $\tilde Z$ for $3\le i\le n+1$, we have
\begin{equation}\label{bs1}
{\text{Bs}}\,|\tilde{W}_m|\subset \bigcup_{3\le i\le n+1}(C_i\cup\ol{C}_i).
\end{equation} 
 
 The restriction of the meromorphic map $\tilde{Z}\to\mathscr T$  onto $E_1$ (and $\ol{E}_1$) is a birational morphism
 (since the composition $Z_1\to Z\to\mathscr T$ is already birational on $E_1$).
 Since $E_1$ (and $\ol{E}_1$) is a  smooth surface, this restriction is a resolution of the minitwistor space $\mathscr T$.
Although $E_1$ is not a minimal surface, the restriction can be seen to be a minimal resolution of $\mathscr T$.
In fact, the restriction of the blowing-up $\mu_2:Z_2\to Z_1$ resolves the conjugate pair of singularities $P_{\infty}$ and $\ol{P}_{\infty}$ ((a) in Proposition \ref{prop-sing_T}).
The remaining blowing-ups $\tilde Z\to\cdots\cdots\to Z_2$ resolve
the remaining real $A_{l_i-1}$-singularities ((b) of Proposition \ref{prop-sing_T}),
which has to be minimal since by our choice, precisely $(l_i-1)$ rational curves are inserted through our blow-ups, for each $i$ satisfying $l_i\ge 2$.
Thus, if $\tilde{\mathscr T}\to \mathscr T$ denotes the minimal resolution of $\mathscr T$, we obtain the following commutative diagram of the meromorphic maps
\begin{equation}\label{cd3}
 \CD
\tilde{Z}@>>>Z\\
 @V\tilde{\Phi}_m^{G_1}VV @VV{\Phi}_m^{G_1}V\\
\tilde{\mathscr T}@>>>\mathscr T\\
 \endCD
 \end{equation}
where the horizontal arrows are the above sequence of blowing-ups and the minimal resolution respectively, and  the map $\tilde{\Phi}_m^{G_1}$ is a meromorphic map which is uniquely determined by the commutativity of the diagram.
Since the restrictions of the composition $\tilde Z\to Z\to\mathscr T$ onto $E_1$ and $\ol{E}_1$ induce the resolution of $\mathscr T$ as above, 
the map $\tilde{\Phi}_m^{G_1}$ also has no indeterminacy locus in a neighborhood of $E_1\cup\ol{E}_1$.

On the other hand, as for the normal bundle of $E_1$ in $\tilde{Z}$, we note that at first its normal bundle in $Z_1$ is isomorphic to $\mathscr O(l+1,-1)$, where $l=(C_1)^2_S$ as before, and $\mathscr O(0,1)$ is a fiber class of the projection $E_1\to C_1$.
Then since all blowing-ups for obtaining $\tilde{Z}$ are performed along curves which intersect $E_1\cup\ol{E}_1$ at  points, we have
\begin{equation}\label{nb1}
N_{E_1/\tilde{Z}}\simeq \nu^*\mathscr O(l+1,-1),
\end{equation}
where  $\nu$ is the restriction of the composition of blowing-ups $\tilde{Z}\to\cdots\to Z_2\to Z_1$ onto $E_1$.

\subsection{Projective models of the twistor spaces of Joyce metrics as conic bundles}
We recall that the sequence of blowing-ups $\tilde{Z}\to\cdots\to Z_1\to Z$ given in the previous subsection is a partial elimination of the indeterminacy locus of the meromorphic map $\Phi_m^{G_1}$ associated to the system $|W_m|$ on the twistor space $Z$,
and we are writing $|\tilde{W}_m|$ to mean the linear system on $\tilde{Z}$ corresponding to the original system $|W_m|$ on $Z$.
By \eqref{bs1}, the indeterminacy locus of the map $\tilde{\Phi}_m^{G_1}:\tilde{Z}\to\tilde{\mathscr T}$ is contained in 
 $\cup_{3\le i\le n+1}(C_i\cup\ol C_i)$, where we are viewing these curves as those in $\tilde{Z}$ as before.
Let $\tilde{Z}'\to \tilde{Z}$ be the composition of a sequence of blowing-ups which eliminates all these indeterminacy locus.
We may suppose that the exceptional divisors of these blow-ups are over the curves $C_i$ and $\ol{C}_i$ with $3\le i\le n+1$.
In particular, they are disjoint from the divisors $E_1$ and $\ol E_1$.
Further we can also suppose that all  centers of the blow-ups are $G_1$-invariant.
Then $\tilde Z'\to \tilde Z$ is biholomorphic in a neighborhood of $E_1\cup \ol E_1$.
 Thus  by taking a composition with $\tilde Z\to\tilde{\mathscr T}$, we obtain a holomorphic map
\begin{equation}
\tilde{\Phi}':\tilde{Z}'\to \tilde{\mathscr T},
\end{equation}
which is bimeromorphic to the original meromorphic map $\Phi^{G_1}_m:Z\to\mathscr T$.

In order to obtain a projective model of the twistor space $Z$ as a conic bundle over the resolved minitwistor space $\tilde{\mathscr T}$, we consider the short exact sequence
\begin{align}\label{ses1}
0\,\lra\, \mathscr O_{\tilde Z'}\,\lra\,\mathscr O_{\tilde Z'}(E_1+\ol{E}_1)\,\lra\,N_{E_1/\tilde{Z'}}\oplus N_{\ol{E}_1/\tilde{Z'}}\,\lra\,0.
\end{align}
Since any fiber of $\tilde{\Phi}'$ is at worst a string of rational curves whose end components intersect $E_1$ and $\ol{E}_1$ respectively, a direct image sheaf satisfies $R^1\tilde{\Phi}'_*\mathscr O_{\tilde Z'}=0$.
Hence by taking the direct image of the sequence \eqref{ses1} and recalling that $E_1$ and $\ol{E}_1$ are sections of $\tilde{\Phi}'$, we obtain the short exact sequence
\begin{align}\label{ses2}
0\,\lra\, \mathscr O_{\tilde{\mathscr  T}}\,\lra\,\tilde\Phi'_*\mathscr O_{\tilde Z'}(E_1+\ol{E}_1)\,\lra\,N_{E_1/\tilde{Z'}}\oplus N_{\ol{E}_1/\tilde{Z}'}\,\lra\,0,
\end{align}
where we regard $N_{E_1/\tilde{Z'}}$ and $N_{\ol{E}_1/\tilde{Z}'}$ as line bundles over $\tilde{\mathscr T}$ by using the isomorphisms $\tilde{\Phi}'|_{E_1}$ and $\tilde{\Phi}'|_{\ol E_1}$.
Since $\tilde Z'$ and $\tilde Z$ are isomorphic in a neighborhood of $E_1\cup\ol{E}_1$, we have $N_{E_1/\tilde Z'}\simeq N_{E_1/\tilde Z}$ and $N_{\ol E_1/\tilde Z'}\simeq N_{\ol E_1/\tilde Z}$.
The same is true for the normal bundle of $\ol{E}_1$.
Therefore  by \eqref{nb1}  we have $H^1(N_{E_1/\tilde{Z'}})=H^1(N_{\ol{E}_1/\tilde{Z'}})=0$
as a consequence of the fact that the second entry of \eqref{nb1} is $-1$.
This means that  the sequence \eqref{ses2} splits.
Hence we obtain
\begin{align}\label{dis1}
\tilde\Phi'_*\mathscr O(E_1+\ol{E}_1)
\simeq
N_{E_1/\tilde{Z}'}\oplus N_{\ol{E}_1/\tilde{Z}'}\oplus\mathscr O_{\tilde{\mathscr T}}.
\end{align}
Let
\begin{equation}\label{p2bdle}
\mu: \tilde{Z}'\lra\mathbf{P}( N_{E_1/\tilde{Z}'}^{\vee}\oplus N_{\ol{E}_1/\tilde Z'}^{\vee}\oplus\mathscr O_{\tilde{\mathscr T}})
\end{equation}
be the relative meromorphic map over $\tilde{\mathscr T}$ associated to the pair of the morphism $\tilde\Phi'$ and the line bundle $\mathscr O(E_1+\ol{E}_1)$. 
Since the restriction of $\mathscr O(E_1+\ol{E}_1)$ to smooth fibers are isomorphic to $\mathscr O(2)$, the map $\mu$ is a bimeromorphic map whose image is a conic bundle in the $\mathbf{CP}^2$-bundle of \eqref{p2bdle}.
We denote this conic bundle by $p:X\to\tilde{\mathscr T}$.
The discriminant locus of $p$ is a member of the linear system $|N_{E_1/\tilde{Z}'}^{\vee}\otimes N_{\ol{E}_1/\tilde Z'}^{\vee}|$.
Let $\tilde{\pi}_m:\tilde{\mathscr T}\to\Lambda_m$ be the composition morphism of the minimal resolution $\tilde{\mathscr T}\to\mathscr T$ and $\pi_m:\mathscr T\to\Lambda_m$.
$\tilde{\pi}_m$ factors as $\tilde{\mathscr T}\to\hat{\mathscr T}\to\Lambda_m$, where $\hat{\mathscr T}\to\Lambda_m$ is the composition 
$\hat{\mathscr T}\to\mathscr T\to\Lambda_m$.
 Since $\hat{\mathscr T}\to\Lambda_m$  has reducible fibers precisely over $(1,\lambda_i,\cdots,\lambda_i^m)\in\Lambda_m$ with $l_i>0$ by Proposition \ref{prop-rf1}, and since $\tilde{\mathscr T}\to\hat{\mathscr T}$ gives the minimal resolution of the $A_{l_i-1}$-singularities of $\hat{\mathscr T}$, $\tilde{\pi}_m$ has reducible fibers precisely over the same points and they consist of $2+(l_i-1)=(l_i+1)$ rational curves.

\begin{prop}\label{prop-dl1}
Let $p:X\to\tilde{\mathscr T}$ be the conic bundle whose total space is a bimeromorphic projective model of the twistor space of a  Joyce metric as above.
Then its discriminant locus consists of the following.
(a) The two distinguished sections $\Gamma$ and $\ol{\Gamma}$ of the morphism $\tilde{\pi}_m:\tilde{\mathscr T}\to\Lambda_m$ (cf.\,\S 2).
(b) The reducible fibers $\tilde{\pi}_m^{-1}(1,\lambda_i,\lambda_i^2,\cdots,\lambda_i^m)$, where $i$ satisfies $1<i<n+2$ and $l_i>0$.
(c) The irreducible fibers $\tilde{\pi}_m^{-1}(1,\lambda_i,\lambda_i^2,\cdots,\lambda_i^m)$ where $i$ satisfies $1<i<n+2$ and $l_i=0$.
\end{prop}

We note the proposition means that the two reducible fibers over the points $(1,0,0,\cdots,0)$ and $(0,0,\cdots,0,1)\in \Lambda_m$ corresponding to $S_1^++S_1^-$ and $S_{n+2}^++S_{n+2}^-$ are not discriminant curves.

\noindent
{\bf Proof of Proposition \ref{prop-dl1}.}
Since the morphism $\tilde{\Phi}':\tilde{Z}'\to\tilde{\mathscr T}$  cannot have 2-dimensional fibers along some curve on $\tilde{\mathscr T}$, $\tilde{\Phi}'$ is flat outside at most finite number of points on $\tilde{\mathscr T}$.
Hence the restriction of $\mu$ onto a fiber $ (\tilde{\Phi}')^{-1}(y)$ is precisely the meromorphic map associated to the linear system $|\mathscr O(E_1+\ol{E}_1)|_{ (\tilde{\Phi}')^{-1}(y)}|$, except at most finite number of points $y\in\tilde{\mathscr T}$.
In particular, if a divisor $E$ on $\tilde{Z}'$  is an irreducible component of a divisor $ (\tilde{\Phi}')^{-1}(D)$ for some curve $D\subset\tilde{\mathscr T}$, and if  $\tilde{\Phi}'(E)=D$, then 
$\mu(E)$ remains  a divisor on $X$.
In the following argument, we frequently use this fact.

By our choice of the blow-up sequence $\tilde{Z}\to\cdots \to Z_1\to Z$, 
either $\tilde{\Phi}^{G_1}_m(E_2)=\tilde{\Phi}^{G_1}_m(E_{n+2})=\Gamma$ or
$\tilde{\Phi}^{G_1}_m(E_2)=\tilde{\Phi}^{G_1}_m(E_{n+2})=\ol{\Gamma}$ holds.
By a possible renaming for $\Gamma$ and $\ol{\Gamma}$, we can suppose the former holds.
Then $\tilde{\Phi}'(E_2)=\tilde{\Phi}'(E_{n+2})=\Gamma$ holds.
Further,  $\tilde{\Phi}'$ gives isomorphism between $E_2\cap E_1\,(\subset \tilde Z')$ and $\Gamma$.
Hence by the above proven facts, the image $\mu(E_2)$ remains a divisor on $X$
satisfying $p(\mu(E_2))=\Gamma$.
Similarly, since  $\tilde{\Phi}'$ gives an isomorphism $E_{n+2}\cap \ol{E}_1\simeq \Gamma$, the image $\mu(E_{n+2})$ remains a divisor
satisfying $p(\mu(\ol E_2))=\ol\Gamma$.
These mean $\mu(E_2\cup E_{n+2})=p^{-1}(\Gamma)$ and that $\mu(E_1)$ and $\mu( E_{n+2})$ are line subbundles of $p$ over $\Gamma$.
Hence $\Gamma$, and therefore also $\ol{\Gamma}$ by reality, are discriminant curves of $p$. 

Next suppose that  $i$ satisfies $1<i<n+2$ and $l_i=1$.
Then by Proposition \ref{prop-rf1}, the fiber $\tilde{\pi}_m^{-1}(1,\lambda_i,\lambda_i^2,\cdots,\lambda_i^m)$ consists of two irreducible components.
Let  $f_i$ and $\ol{f}_i$ be its  irreducible components.
As in the previous subsection let $F_i$ and $\ol{F}_i\subset \tilde Z'$ be the exceptional divisors appeared when we took the blow-up $Z_3\to Z_2$.
Then again by our choice of the blow-up sequence, either
$\tilde{\Phi}'(S_i^-)=\tilde{\Phi}'(\ol{F}_i)=f_i$  or $\tilde{\Phi}'(S_i^-)=\tilde{\Phi}'(\ol{F}_i)=\ol{f}_i$  holds. 
By possible renaming, we can suppose that the former holds.
Then since $\tilde{\Phi}'$ gives isomorphisms $S_i^-\cap E_1\simeq f_i$  and $\ol{F}_i\cap E_1\simeq f_i$, we deduce (again by using the above facts) that 
the image $\mu(S_i^-)$ and $\mu(\ol{F}_i)$ are (mutually distinct) line subbundles of $p$ over $f_i$.
These mean $\mu(S_i^-\cup \ol{F}_i)=p^{-1}(f_i)$.
Hence $f_i$, and therefore also $\ol{f}_i$ by reality, are discriminant curves of $p$. 

Next suppose that  $i$ satisfies $1<i<n+2$ and $l_i>1$.
Then the fiber $\tilde{\pi}_m^{-1}(1,\lambda_i,\lambda_i^2,\cdots,\lambda_i^m)$ is a chain of $(l_i+1)$ rational curves.
Recall that  the two end components  intersect $\Gamma$ and $\ol{\Gamma}$ respectively, and that the other intermediate components are the exceptional curves of the minimal resolution of $A_{l_i-1}$-singularity of the original minitwistor space $\mathscr T$.
Again by our way of blowing-ups, the images $\tilde{\Phi}'(S_i^-)$ and $\tilde{\Phi}'(\ol{F}_i)$ 
are same end component of the chain.
Also $S_i^-\cap E_1$ and $\ol{F}_i\cap E_1$ are mapped by $\tilde{\Phi}'$ isomorphically onto this end component.
Hence we  deduce that this end component of the chain is actually a discriminant curve of $p$.
By reality, the other (conjugate) end component is also a discriminant curve of $p$.
On the other hand,  all the remaining intermediate components of the chain are the images of the two exceptional divisors appearing when we performed the blowing-up $\tilde{Z}\to\cdots\to Z_4\to Z_3$ in the previous subsection, and all of them  intersect $E_1$ or $\ol{E}_1$.
These mean that the intermediate components of the chain are also discriminant curves.
Thus we obtain that all  curves in  (b) in the proposition are discriminant curves.

Next suppose that $i$ satisfies $1<i<n+2$ and $l_i=0$.
Then the fiber $\pi_m^{-1}(1,\lambda_i,\lambda_i^2,\cdots,\lambda_i^m)$ is irreducible
 by Proposition \ref{prop-rf1} and is exactly $\Phi_m^{G_1}(L_i)$ by the proof of Proposition \ref{prop-irf}.
Obviously we have $\tilde{\Phi}'(L_i)=\Phi_m^{G_1}(L_i)$ by $\tilde{\mathscr T}\to\mathscr T$.
Again since $\tilde{\Phi}'$ gives isomorphisms $S_i^+\cap E_1\simeq \tilde{\Phi}'(L_i)$ and $S_i^-\cap E_1\simeq \tilde{\Phi}'(L_i)$, 
the images $\mu(S_i^+)$ and $\mu(S_i^-)$ are mutually different line subbundles of $p$ over the fiber  $\tilde{\Phi}'(L_i)$.
These mean that the irreducible fiber in the item (c) is also a discriminant curve.

To complete a proof, it remains to show that there exists no discriminant curve other than those in (a), (b) and (c).
Let $D$ be such an irreducible curve.
Then since the conic bundle map $Y\to\tilde{\mathscr T}$ is $G$-equivariant, $D$ must be $G$-invariant.
Therefore $D$ is an irreducible component of a fiber of the projection $\tilde{\pi}_m:\tilde{\mathscr T}\to\Lambda_m$.
Suppose that this fiber is irreducible.
Then the inverse image $p^{-1}(D)$  consists of two irreducible components (since $p$ is a conic bundle).
Consider the strict transforms of these two divisors into the original twistor space $Z$.
Then since $D$ is a fiber of $\tilde{\mathscr T}\to\Lambda_m$, the sum of these two transforms is a member of the pencil $|F|^G$.
Since we are supposing that $D$ is not in the item (b) and (c), this means that there exists a reducible member of the pencil $|F|^G$ other than $S_i^++S_i^-$, $1\le i\le n+2$.
Such a divisor does not exist and hence we deduce that $D$ cannot be a reducible fiber.
Suppose finally that $D$ is contained in a reducible fiber of  $\tilde{\mathscr T}\to\Lambda_m$.
Then by assumption that $D$ is not (b) nor (c),
the reducible fiber must be over the point $(1,0,0,\cdots,0)$ or $(0,0,\cdots,0,1)\in \Lambda_m$.
Since $l_1=l_{n+2}=1$, these fibers consist of two irreducible components.
We consider the reducible fiber over $(1,0,0,\cdots,0)$.
Then the divisor $S_1^+$ and $S_1^-$ in $\tilde Z'$ are mapped surjectively to mutually different irreducible components of this fiber by the morphism $\tilde Z'\to\tilde{\mathscr T}$ since it is true for the meromorphic map $Z\to \mathscr T$.
By our choice of the blowing-ups, both $S_1^+$ and $S_1^-$ intersect $E_1$ and $\ol{E}_1$ along  curves which are mapped biholomorphically to the irreducible components of the fiber.
This means that the restriction of $\mu:\tilde Z'\to X$ to $S_1^+$ and $S_1^-$ are birational and the images are exactly the inverse images of the irreducible components under the map $X\to\tilde{\mathscr T}$.
Hence the fiber over $(1,0,\cdots,0)$ cannot be discriminant curves of $\mu$.
The same is also true for the fiber over the point $(0,0,\cdots,0,1)$ as one can see by replacing $S_1^+$ and $S_1^-$ by $S_{n+2}^+$ and $S_{n+2}^-$ respectively in the above argument.
\proofend

\vspace{2mm}
As is already mentioned, the discriminant locus of the conic bundle $p:X\to\tilde{\mathscr T}$ is a member of the system $|N_{E_1/\tilde{Z}'}^{\vee}\otimes N_{\ol{E}_1/\tilde Z'}^{\vee}|$, where the line bundles $N_{E_1/\tilde{Z}'}^{\vee}$ and $ N_{\ol{E}_1/\tilde Z'}^{\vee}$ can be explicitly described thanks to \eqref{nb1}.
However, a computation shows that the sum of all discriminant curves obtained in Proposition \ref{prop-dl1} is not a member of the system in general;
usually the subtraction of the sum from the system has an effective member.
This means that the discriminant curves are non-reduced in general.

 \section{Equivariant deformations of the twistor spaces of Joyce metrics}
 \subsection{}
 In this subsection we investigate deformations of the twistor spaces of Joyce metrics which preserve a $G_1$-action, where $G_1$ is a prescribed $\mathbf C^*$-subgroup of $G$ as in Sections 2 and 3.
We will give a sufficient condition for the subgroup $G_1$ for which the twistor space admits a $G_1$-equivariant but non $G$-equivariant deformation.
Since such deformations are well understood for $G_1$ acting semi-freely (on LeBrun twistor spaces) (\cite{LB91,LB93}), we suppose that $G_1$ does not act semi-freely.
(But we do not suppose that the twistor spaces are not LeBrun's one.)

As in the previous sections,  $Z$ denotes the twistor space of a Joyce metric on $n\mathbf{CP}^2$, and $C$ denotes the cycle of $2(n+2)$ rational curves which is the unique $G$-invariant anticanonical curve of a smooth member $S$ of the pencil $|F|^G$.
To state our result precisely, we introduce the following.

\begin{definition}\label{def-reg}{\em
Let $S$ and $C$ be as above, and $G_i$ the $\mathbf C^*$-subgroup of $G$
which fixes $C_i$ and $\ol{C}_i$.
An irreducible component $C_j$ (or equivalently, $\ol{C}_j$) ($j\neq i$) is said to be {\em  regular} with respect to $G_i$ if the isotropy subgroup of $G_i$-action on $C_j$ is the  identity only.
We say that the component is {\em irregular} with respect to $G_i$ if it is not regular.
}
\end{definition}

\noindent
We simply say that a component $C_j$ (or equivalently $\ol{C}_j$) is regular if the subgroup $G_i$ or the component $C_i$  is obvious from the context.

In the following, we again adopt the numbering for the irreducible components of  $C$ for which the chosen component is $C_1$, so that the $\mathbf C^*$-subgroup of $G$ is $G_1$.
If still denoting the two reducible fibers of the $G_1$-quotient map $S\to\mathbf{CP}^1$ (of \eqref{quot_S}) by $f$ and $\ol f$ and writing them as in \eqref{f_1}, then $k_i\in\mathbf Z_{>0}$ coincides with the order of the isotropy subgroup of
$G_1$-action on $C_i$.
Hence $C_i$ is regular  (with respect to $G_1)$ iff $k_i=1$.
Therefore the two components $C_2$ and $\ol{C}_{n+2}$
are always regular.
If we realize the toric surface $S$ as a succession of $G$-equivariant blow-ups of $\mathbf{CP}^1\times\mathbf{CP}^1$ preserving the real structure and  satisfying the condition that the component $C_1$ (and hence $\ol{C}_1$ also) is not an exceptional curve,
then irregular components appears only when we blow-up an isolated fixed point of $G_1$-action.
In other words, the exceptional curve of a blow-up becomes regular only when
the blown-up point is on $C_1\cup\ol C_1$.
Since we are supposing non semi-freeness for $G_1$, there always exists an irregular component $C_j$, $3\le j\le n+1$.

Next by using this regularity, we define two integers $r$ and $s$ for an arbitrary pair $(Z,G_1)$.
As noted above, the two components $C_2$ and $C_{n+2}$ are regular.
Let $r$ be the maximum integer satisfying the condition that $C_j$ is regular   for all $2\le j\le r$.
Let $s$ be the minimum integer satisfying the condition that $C_j$ is regular  for all $s\le j\le n+2$.
Since we have supposed non semi-freeness for $G_1$, we have $(2
\le )\,r<s\,(\le n+2)$.
The unions 
$C_2\cup C_3\cup\cdots \cup C_{r}$ and $C_{s}\cup\cdots\cup C_{n+2}$ are connected components of the union of all regular components, which intersect $C_1$ and $\ol{C}_1$ respectively.
The two integers $r$ and $s$ are uniquely determined once the component $C_1$ is given, and  easily computable since the weights $k_i$  are easily computable.
(We note that in general, there can exist a regular component $C_j$ satisfying $r<j<s$.)

The following result gives a sufficient condition for  the subgroup $G_1$ to have a $G_1$-equivariant but non-$G$-equivariant deformation of $Z$.
We recall that we have explicitly constructed the minitwistor space $\mathscr T$ that can be regarded as a quotient space of the $G_1$-action (Definition \ref{def-mt}).

 \begin{thm}\label{thm-1}
Let  $Z$ be the twistor space of an arbitrary Joyce metric on $n\mathbf{CP}^2$.
Let $C_1\subset C$ be any irreducible component of $C$.
Suppose that the isotropy subgroup $G_1$ of $C_1$ acts non semi-freely, so that $n\ge 2$.
Let $\mathscr T$ be the minitwistor space of $Z$ with respect to $G_1$.
Then if the two integers $r$ and $s$ defined above satisfy an inequality
\begin{equation}\label{ineq-key}
n+r-s>0,
\end{equation}
then $Z$ can be $G_1$-equivariantly deformed into a twistor space $Z_b$ satisfying the following.
(a) $Z_b$ does not admit an effective $G$-action.
(b) The linear system $|mF|^{G_1}$ on $Z_b$ remains $(m+2)$-dimensional, and the image of the associated map is still  $\mathscr T$. 
 \end{thm}
 
We note that since $2\le r<s\le n+2$, we always have $n+r-s\ge 0$.
 We also note that the inequality \eqref{ineq-key} is a condition for the diffeomorphism type of $U(1)$-action on $n\mathbf{CP}^2$ corresponding to $G_1$-action.
That is, if  $K_1\,(\simeq U(1))$ is the subgroup of $K$ corresponding to $G_1$, $(n+r-s+2)$ is the number of $K$-invariant spheres in $n\mathbf{CP}^2$  which are acted by weight one by $K_1$ and which are joined to the $K_1$-fixed sphere through the spheres.

%
 
 We prove Theorem \ref{thm-1} in the following manner.

 \begin{enumerate}
 
 \item[$1^{\circ}$]
 We show that, under the condition  \eqref{ineq-key},
 we explicitly construct a $G_1$-equivariant deformation of the toric surface $S\in |F|^G$
 which breaks the structure of toric surface, and for which the components $C_1$  and $\ol{C}_1$ survive.
\item[$2^{\circ}$] We consider deformations of the pair $(Z,S)$ and apply an equivariant version of a theorem of E.\,Horikawa \cite{Hor76} to conclude that there exists a $G_1$-equivariant deformation of the pair such that, if restricted to $S$, the deformation coincides with  the $G_1$-equivariant deformation constructed in (1).
By looking the structure of the non-toric surface, we see that the deformed twistor space cannot admit an effective $G$-action.
\item[$3^{\circ}$]
Let $m$ be the integer in Definition \ref{def-m}, so that the linear system $|W_m|$ on the twistor space $Z$ (of a Joyce metric) induces a meromorphic quotient map $\Phi_m^{G_1}:Z\to\mathscr T\subset\mathbf{CP}^{m+2}$.
By a similar argument we  employed to $Z$, we show that the system $|mF|^{G_1}$ on the deformed twistor space   is still  $(m+2)$-dimensional whose  image  of the associated map is  a complex surface whose defining equation is the same as the original minitwistor space.
 \end{enumerate}
  
\noindent
{\bf{Proof of Theorem \ref{thm-1}.}}
$1^{\circ}$ By definition of the integers $r$ and $s$, a component $C_j$ is regular with respect to $G_1$ as long as $2\le j\le r$ and $s\le j\le n+2$.
So we have $(r-1)+(n-s+3)=(n+r-s+2)$ regular components at least.
By \eqref{ineq-key}, we have $n+r-s+2>2$.
This means that either  $C_3$ or $C_{n+1}$ is a regular component at least.

Suppose that $C_3$ is a regular component (so that $r\ge3$).
We show that $C_2,C_3,\cdots,C_{r-1}$ can be successively blown-down in $S$ in this order.
For this, we fix any sequence of blowing-downs of $S$ to a minimal surface which does not contract $C_1$ nor $\ol{C}_1$ and view the sequence as a sequence of blowing-ups.
Then an irreducible component $C_j$ ($j\neq 1$) is regular iff it is the exceptional curve of a blowing-up at $C_1$ or $\ol{C}_1$, or it is not contracted by the sequence of blowing downs.
If $C_3$ is contracted to a point by these blowing-downs, since $C_3$ is supposed to be regular and there exists irregular component $C_k$ satisfying $k>3$ by assumption, $C_3$ is an exceptional curve which arises when we blow-up a $G$-fixed point of $C_1$.
This means $C_2^2=-1$.
If $C_3$ is not contracted to a point by these blowing-downs, $C_2^2=-1$ holds obviously.
Hence we  obtain $C_2^2=-1$ (if $C_3$ is a regular component).
So we blow-down $C_2$. Then $C_3$ is of course still regular, intersecting $C_1$.
If $C_4$ is also a regular component (i.e. if $r\ge 4$), we employ the same argument by replacing  $C_2$ and $C_3$ by $C_3$ and $C_4$ respectively.
Consequently we obtain $C_3^2=-1$ on the blown-down surface.
So we can blow-down $C_3$.
Repeating this argument, we conclude that $C_2, C_3,\cdots,C_{r-1}$ can be successively blown-down in this order.
By reality, the same is true for $\ol{C}_2,\ol{C}_3,\cdots,\ol{C}_{r-1}$.
The same argument shows that if $C_{n+1}$ is a regular component (so that $s\le n+1$), then
$C_{n+2}, C_{n+1}, \cdots, C_{s+1}$ and also 
$\ol{C}_{n+2}, \ol{C}_{n+1}, \cdots, \ol{C}_{s+1}$ can be successively blown-down in this order.
Summing all these up, we can successively blow-down $2(n+r-s)$ components.
Let $S\to \ol{S}$ be this blowing-down contracting $2(n+r-s)$ curves.
The surface $\ol S$ is uniquely determined from the toric surface $S$ and the subgroup $G_1$.
(The condition \eqref{ineq-key} guarantees that $S\to\ol{S}$ is a non-trivial map.)
The curves $C_2, C_3,\cdots,C_{r-1}$ are contracted to a $G$-fixed point of $C_1$, and 
$\ol{C}_{n+2},\ol{C}_{n+1},\cdots,\ol{C}_{s+1}$ are mapped to another $G$-fixed point on $C_1$.
Similarly, the curves $\ol{C}_2,\ol{C}_3,\cdots,\ol{C}_{r-1}$ are contracted to a $G$-fixed point of $\ol{C}_1$, and $C_{n+2},C_{n+1},\cdots,C_{s+1}$ are contracted to another $G$-fixed point on $\ol{C}_1$.

Next by using this blowing-down $S\to\ol{S}$ 
we concretely construct a $G_1$-equivariant deformation of the  surface $S$.
For this purpose,
put $B:=C_1^{\times (n+r-s)}\times \ol C_1^{\times (n+r-s)}$ and consider a direct product $\ol{S}\times B$.
We view this as a trivial deformation of $\ol S$.
This trivial $\ol S$-bundle has a natural (product) $G$-action, since $\ol{S}$ is  a toric surface and $C_1$ and $\ol{C}_1$ are $G$-invariant.
The subgroup $G_1$ acts trivially on the base space $B$.
Further, this bundle has a tautological multi-section (which is generically $2(n+r-s)$ to $1$) determined from the inclusions $C_1\subset \ol S$ and $\ol C_1\subset \ol S$.
Let $\mathscr S\to \ol S\times B$ be the blowing-up along this multi-section and 
$p:\mathscr S\to B$ the composition of $\mathscr S\to \ol S\times B$ and $\ol S\times B\to B$.
Since the multi-section is $G$-invariant, $p:\mathscr S\to B$ has a natural $G$-action, where $G_1$ acts trivially on $B$.
If $P, Q\in C_1\subset \ol S$ are the images of $C_2\cup\cdots \cup C_{r-1}$ and $\ol{C}_{n+2}\cup\cdots\cup\ol C_{s+1}$ under the blowing-down $S\to\ol S$ respectively, 
the point
$b_0:=(P,\cdots, P,Q,\cdots, Q;\ol P,\cdots,\ol P,\ol Q,\cdots,\ol Q)\in B$
satisfies $p^{-1}(b_0)\simeq S$.
Thus  $p$ is a $G_1$-equivariant deformation of $S$.
This is the required  deformation.
Since  the $G_1$-quotient map $\ol S\to\mathbf{CP}^1$ has at least 2 reducible fibers, the $G_1$-quotient map $S_b\to\mathbf{CP}^1$ has at least 4 reducible fibers  for general $b\in B$.
Hence $S_b$ is not  a toric surface for general $b\in B$.
We also note that since $C_i$ and $\ol C_i$ are not contracted by $S\to \ol S$ for $i=1$ and $r\le i\le s$,  these curves still make sense on $S_b$ for any $b\in B$.
Then it can be readily verified that  the curve
\begin{equation}\label{ac3}
C_1+\sum_{r\le i\le s}C_i+\ol{C}_1+\sum_{r\le i\le s}\ol{C}_i
\end{equation}
is an anticanonical curve on $S_b$ for general $b\in B$, forming the cycle (like $C$ in $S$).


$2^{\circ}$ Now since the the twistor space $Z$  is Moishezon, we have the $H^2(Z,\Theta_Z\otimes\mathscr O(-S))=0$ (\cite[Lemma 1.9]{C91}).
Then by a theorem of Horikawa \cite{Hor76}, this implies so called the co-stability of $S$; namely for any deformation $q:\mathscr S\to T$ of $S=q^{-1}(t_0)$ with $t_0\in T$, there exist a deformation $\tilde{q}:\mathscr Z\to T'\subset T$ with $T'$ open containing $t_0$ with $\tilde{q}^{-1}(t_0)\simeq Z$, and an inclusion $i:\mathscr S\subset\mathscr Z$ over $T'$, which satisfy $q=\tilde q\circ i$.
Applying this to the present deformation $p:\mathscr S\to B$,
we obtain a deformation $\tilde{p}:\mathscr Z \to B'\subset B$ with $B'$ open containing $b_0$ with $\tilde p^{-1}(b_0)\simeq Z$, as well as an inclusion $i:\mathscr S\subset \mathscr Z$ over $B'$ satisfying ${p}=\tilde p\circ i$.
Further, since $(Z,S)$ admits $G_1$-action and $p$ is a $G_1$-equivariant deformation, 
$\tilde p$ can be supposed to be a $G_1$-equivariant deformation of $Z$.
If we restrict $\tilde p$ to a real locus $(B')^{\sigma}$ of $B'$ (which is necessarily a smooth submanifold of the half-dimension of $B'$), we obtain a deformation of $Z$ which preserves $G_1$-action, the surface $S$ and the real structure.
Further, since any small deformation of a compact twistor space remains a twistor space, we can suppose that all fibers over $(B')^{\sigma}$ are twistor spaces, after possible shrinking of $(B')^{\sigma}$.

Take a point $b\in(B')^{\sigma}$ such that $S_b$ is not a toric surface.
Then $Z_b=\tilde{p}^{-1}(b)$ is a twistor space  equipped with a $G_1$-action and a $G_1$-invariant surface $S_b=p^{-1}(b)\subset Z_b$.
Since the Picard group of the twistor spaces is isomorphic to the (integral) second cohomology group, $S_b$ is still a member of $|F|$ on $Z_b$.
We show that $Z_b$ does not admit an effective $G$-action.
If $Z_b$ admits such an action, $Z_b$ must be a LeBrun twistor spaces since any smooth member of $|F|$ in non-LeBrun twistor spaces of Joyce metrics is a toric surface.
Assume that the original space $Z$ is not a LeBrun twistor space.
Then since there exists no component $C_i$ of $C$ satisfying $(C_i)^2_S=-n$, the same must be true for $S_b$ about the cycle \eqref{ac3}.
On the other hand, for any LeBrun twistor spaces, any anticanonical curve of a smooth member of $|F|$ has an irreducible component $C_i$ satisfying $C_i^2=-n$.
(It is concretely given by the irreducible component of Bs\,$|F|$.)
This is a contradiction and hence $Z_b$ does not admit an effective $G$-action.
Next assume that $Z$ is a LeBrun twistor space.
Then we can verify that among $(n+1)$ subgroups acting non semi-freely, there exist precisely two subgroups satisfying the condition \eqref{ineq-key}
(cf.\,\S 5.2).
Applying the above construction for these subgroups to $S\in |F|$, we obtain a non-toric surface $S_b$.
It can be seen that the cycle  has no $(-n)$-component (though $(-n+1)$-component exists).
This is again a contradiction and we conclude that $Z_b$ does not admit an effective $G$-action.

$3^{\circ}$
Let  $m$ be the integer in Definition \ref{def-m} applied to $(S,G_1)$.
As is already proved, the system $|W_m|$ on  $Z$ induces the surjective meromorphic map 
$\Phi_m^{G_1}:Z\to\mathscr T$ which can be regarded as a quotient map of the $G_1$-action on $Z$.
Now using the explicit structure of the surface $S_b\in |F|$ on $Z_b$, we show that the system $|mF|$ on $Z_b$ still has $(m+2)$-dimensional subsystem which induces a meromorphic map that can yet be regarded as  a quotient map for the $G_1$-action on $Z_b$.
To see this, as before, let $f$ and $\ol f$ be the reducible fibers \eqref{f_1} of the $G_1$-quotient map $S\to\mathbf{CP}^1$ for the toric surface.
As is already noted,  the curves $C_j$ and $\ol C_j$ naturally make sense on $S_b$ as long as $j=1$ or $r\le j\le s$.
We consider the quotient morphism $S_b\to\mathbf{CP}^1$ for the $G_1$-action on $S_b$.
This has distinguished two reducible fibers which correspond to the two reducible fibers $f$ and ${\ol f}$ of $S\to\mathbf{CP}^1$.
We write them $f_b$ and $\ol{f}_b$ respectively.
We can show that 
\begin{equation}\label{fiber1}
f_b=\sum_{r\le i\le s}k_iC_i
\end{equation}
hold, by using the fact that the cohomology class $[f_b]$ of  fibers of the quotient map $S_b\to\mathbf{CP}^1$ is characterized by the property that and $f_b\cdot C_1=1$ and $f_b\cdot D=0$ , where $D$ is any irreducible components of arbitrary fibers of $S_b\to\mathbf{CP}^1$.
Since this is a claim for the rational surface whose structure is explicitly given, and since the computations is long, we omit the details here.

Next by using \eqref{fiber1} we find a divisor $Y_b\in |mF|$ on $Z_b$ which is again a sum of degree one divisors.
As in the case of the twistor space of a Joyce metric, let $L_i\subset Z_b$ be the twistor line which goes through the intersection point $C_1\cap C_r$ for the case $i=1$, $C_i\cap C_{i+1}$ for the case $r\le i< s$, and $C_{s}\cap \ol{C}_1$ for the case $i=s$.
All these twistor lines are $G_1$-invariant.
Then by \cite[Proposition 3.7]{Kr98}, there exists a set of $2(s-r+2)$ degree one divisors $\{S_i^+, S_i^-\set i=1{\text{ or }} r\le i\le s\}$ on $Z_b$ satisfying $\ol{S}_i^+=S_i^-$ and $S_i^+\cap S_i^-=L_i$ for each $i$.
Since the restriction $(S_i^++S_i^-)|_{S_b}$ is necessarily an anticanonical curve \eqref{ac3}  and since $S_i^+$ and $S_i^-$ are mutually conjugate,
the restrictions $S_i^+|_{S_b}$ and  $S_i^-|_{S_b}$ are precisely two connected halves of the cycle \eqref{ac3} divided by the twistor line $L_i$.
We again make a distinction between $S_i^+$ and $S_i^-$ by imposing that $S_i^+$ contains the component $\ol{C}_1$.  
Also we still  write $C$ for the cycle \eqref{ac3} on $Z_b$.

Now we consider the two curves $mC+f_b-\ol{f}_b$ and $mC-f_b+\ol{f}_b$ on  $S_b$.
Since $f_b$ is explicitly given by \eqref{fiber1} and the coefficients $k_i$  are exactly the same  as those for the reducible fiber $f$ on the toric surface $S$, the two curves $mC+f_b-\ol{f}_b$ and $mC-f_b+\ol{f}_b$ are effective curves (because this is true for the curves $mC+f-{\ol f}$ and $mC-f+{\ol f}$).
Since $f_b$ and $\ol{f}_b$ are linearly equivalent on $S_b$, we have $mC+f_b-\ol{f}_b\in |mK^{-1}_{S_b}|$ and $mC-f_b+\ol{f}_b\in |mK^{-1}_{S_b}|$.
Further, by \eqref{fiber1} we have
\begin{equation}\label{33}
mC+f_b-\ol{f}_b=mC_1+\sum_{r\le i\le s}(m+k_i)C_i+m\ol{C}_1+\sum_{r\le i\le s}(m-k_i)\ol{C}_i.
\end{equation}
Next we apply Procedure (A) to the sequence $(k_r,k_{r+1},\cdots,k_{s})$.
This sequence is of course a part of 
the sequence  $(k_2,k_3,\cdots,k_{n+2})$ for the case of the original twistor space $Z$ of a Joyce metric.
Further, since $k_2=\cdots=k_r=1$ and 
$k_{s}=k_{s+1}=\cdots=k_{n+2}=1$, 
the number $m$ defined in Definition \ref{def-m} are exactly the same for the sequence  $(k_2,k_3,\cdots,k_{n+2})$ and the subsequence $(k_r,k_{r+1},\cdots,k_{s})$.
Furthermore, the two indices $i=i_l$ and $j=j_l$  in the step (1) of the $l$-th ($1\le l\le m$) application of Procedure (A) are also exactly the same for these two sequences, with only exception in the case $l=m$ (i.e. the final application) in that
 $i_m=r$ and $j_m=s$ for the subsequence.
Then just as before we define the divisor on $Z_b$ by
\begin{equation}\label{Y_b}
Y_b=\sum_{1\le l\le m}(S_{i_l-1}^++S_{j_l}^-).
\end{equation}
Thus at least formally the divisor $Y_b$  is obtained from the divisor $Y$ of \eqref{Y} on $Z$ by just replacing $S_{n+2}^-$ by  $S_{s}^-$.
Then by the same computation in the proof of Proposition \ref{prop-mt3}, by using \eqref{33} we obtain $Y_b|_{S_b}=mC+f_b-\ol{f}_b$.
Since $mC+f_b-\ol{f}_b\in|mK_{S_b}^{-1}|$, we obtain $Y_b\in |mF|$ by the same reason for $Y\in |mF|$ on $Z$.
Hence we also have $\ol{Y}_b\in |mF|$.

As in the case of Joyce metrics, 
let $V_m\subset H^0(Z_b,mF)$ be the $(m+1)$-dimensional linear subspace generated by the image of the map
\begin{equation}
H^0(Z_b,F)\times H^0(Z_b,F)\times\cdots\times H^0(Z_b,F)\,\lra\,H^0(Z_b,mF)
\end{equation}
sending $(s_1,s_2,\cdots,s_m)$ to $s_1\otimes s_2\otimes\cdots\otimes s_m$.
Let $|W_m|$ be the linear system on $Z_b$ generated by $|V_m|$ and $Y_b$ and $\ol{Y}_b$.
This is a $(m+2)$-dimensional  subsystem of $|mF|$.
We have $W_m=H^0(Z_b,mF)^{G_1}$ by the same reason for Proposition \ref{prop-char}.
As in the final part of the proof of Proposition \ref{prop-mt5}, we investigate the meromorphic map associated to this system $|W_m|$ on $Z_b$.
By the inclusion $V_m\subset W_m$ we again  have the following commutative diagram of meromorphic maps
\begin{equation}\label{cd4}
\xymatrix{
   Z_b \ar@{->}[r]^{{\Phi_m^{G_1}}}\ar@{->}[d]_{\Psi_m}  & \mathscr T_b  \ar@{->}[dl]^{\pi_m}\\
   \Lambda_m & \\
}
 \end{equation}
 where $\Psi_m$, $\Phi_m^{G_1}$, $\pi_m$ and $\mathscr T_b=\Phi_m^{G_1}(Z_b)$ have the analogous meaning as those in the diagram \eqref{cd2}.
The restriction of $\Phi_m^{G_1}$ onto a fiber $S_b=\Psi_m^{-1}(\lambda)$ is the meromorphic map associated to the linear system on $S_b$ obtained as the restriction of $|W_m|$ on $S_b$.
The latter system is generated by $Y_b|_{S_b}, \ol{Y_b}|_{S_b}$ and $mC$.
By the same computation of the proof of Proposition \ref{prop-mt2}, the movable part of this system on $S_b$ is base point free, two-dimensional, and its induced morphism is exactly the $G_1$-quotient map $S_b\to\mathbf{CP}^1$, where  $\mathbf{CP}^1$ is a conic embedded in  $\mathbf{CP}^2$.
Hence the meromorphic map $\Phi_m^{G_1}:Z_b\to\mathscr T_b$ can be regarded as a $G_1$-quotient map, and  by $\pi_m$, the image $\mathscr T_b$ has a structure of a (rational) conic bundle over $\Lambda_m$.
Furthermore, since the divisor $Y_b$ has the same form as  $Y$ on $Z$, we can repeat the computation  in the proof of Proposition \ref{prop-mt5} for obtaining the equations \eqref{key1}--\eqref{dcp1} to deduce that  the equation of $\mathscr T_b$ is still given by 
\begin{equation}
z_{m+1}z_{m+2}
=c\,(u_1-\lambda_1'u_{n+2})(u_1-\lambda'_2u_{n+2})^{l_2}(u_1-\lambda'_3u_{n+2})^{l_3}\cdots(u_1-\lambda'_{n+1}u_{n+2})^{l_{n+1}}(u_1-\lambda'_{n+2}u_{n+2}).
\label{dcp4}
\end{equation} 
for some real numbers $\lambda'_1,\cdots,\lambda'_{n+2}$ satisfying 
 $\lambda'_1<\lambda'_2<\cdots<\lambda'_{n+2}$ or
 $\lambda'_1>\lambda'_2>\cdots<\lambda'_{n+2}$,
 where $u_1,u_{n+2},z_{m+1},z_{m+2}$ and $l_k$ have the same meaning as in 
\eqref{key1}--\eqref{dcp1}.
 These mean that the image surface is $G_1$-quotient surface whose structure is the same as the minitwistor space $\mathscr T$ for the original $Z$, with a possible difference of the numbers $\lambda_1<\lambda_2<\cdots<\lambda_{n+2}$.
\proofend

\bigskip
If $(Z,G_1)$ does not satisfy the inequality \eqref{ineq-key}, Theorem \ref{thm-1} does not generate a new twistor space with $G_1$-action.
Hence the problem arises as to whether our minitwistor spaces can be obtained as a quotient space of a  twistor space with only $\mathbf C^*$-action.
The answer is positive as the following result shows.

\begin{thm}\label{thm-2}
Suppose $n\ge 2$ and let
$\mathscr T$ be the minitwistor space of arbitrary Joyce metric on $n\mathbf{CP}^2$ with respect to a $\mathbf C^*$-subgroup $G_1$  in the sense of Definition \ref{def-mt}.
Then for any $n'>n$ there exists a Moishezon twistor space $Z'_b$ on $n'\mathbf{CP}^2$ with $\mathbf C^*$-action which satisfies the following.
(a) $Z'_b$ is not isomorphic to the twistor spaces of Joyce metrics,
(b) By the linear system $|mF|$ on $Z'_b$, $Z'_b$ is mapped onto the minitwistor space $\mathscr T$. 
\end{thm}

In other words, the minitwistor spaces constructed in Section 2 appear not only as those associated to the twistor spaces of Joyce metrics but also as those associated to twistor spaces whose automorphism group is  just $\mathbf C^*$.

\bigskip
\noindent
{\bf Proof.}
Let $(Z,G_1)$ be a pair of the twistor space of a Joyce metric on $n\mathbf{CP}^2$ and the $\mathbf C^*$-subgroup whose minitwistor space is $\mathscr T$.
Let $S\in |F|^G$ be a smooth member, and $C_1$ and $\ol C_1\subset S$ the pair of rational curves fixed by  $G_1$.
Let $r$ and $s$ be the integers appearing in the inequality \eqref{ineq-key}.
Then as is already noted, we have $n+r-s\ge 0$.
For the given integer $n'>n$, let $S'\to S$ be a succession of any $G$-equivariant blow-ups preserving the real structure whose center is always on $C_1\cup\ol C_1$, where the number of times for the blow-ups is precisely $2(n'-n)$.
Then the integer $m'$ obtained by applying Procedure (A) for $(S',G_1)$ is the same as $m$ for $(S,G_1)$, since  the sequence $(k'_2,k'_3,\cdots,k'_{n'+2})$ for $S'$ is obtained from the sequence $(k_2,k_3,\cdots,k_{n+2})$ for $S$ by adding $(n'-n)$ ones for both sides in total.
In particular, if $Z'$ denotes the the twistor space of a Joyce metric on $n'\mathbf{CP}^2$ having $S'$ as a member of the pencil $|F|^G$ on $Z'$, the members $Y'$ and $\ol Y'$ of $|W_m|$ (given by \eqref{Y}) on $Z'$ is of the same form as those of $|W_m|$ on $Z$.
This means that the minitwistor space associated to the pair $(Z', G_1)$ is the same as the one for the pair $(Z,G_1)$.
Further, by our choice of $S'\to S$,  we have $n'+s-s'=n'-n$, which is positive by assumption.
Hence by applying Theorem \ref{thm-1} to $(Z',G_1)$  we obtain the desired twistor space $Z'_b$.
\proofend

 \subsection{Discriminant curves as hyperplane sections of the minitwistor spaces}
In this subsection 
by the method we employed for the twistor spaces of Joyce metrics in \S 3.2, we give projective models of the twistor spaces obtained in the previous subsection and investigate its discriminant curves.
We will show that a `principal part' of the discriminant locus must be hyperplane sections of the minitwistor spaces, with respect to our realization  in $\mathbf{CP}^{m+2}$.

As in the previous section let $|W_m|\,(=H^0(Z_b,mF)^{G_1})$ be a $(m+2)$-dimensional linear system  on  the deformed twistor space $Z_b$.
$|W_m|$ is generated by $Y_b$, $\ol{Y}_b$ and members of $|V_m|$.
Since  the divisors $Y_b$ on $Z_b$ and $Y$ on $Z$ have the same form as in \eqref{Y} and \eqref{Y_b},  in a neighborhood of $C_1\cup\ol{C}_1$, Bs\,$|W_m|$ can be eliminated by the same sequence of blowing-ups we have explicitly given in \S 3.1.
Let $\tilde Z_b\to Z_b$ be this sequence blowing-ups.
Let $|\tilde W_m|$ be the linear system on $\tilde Z_b$ corresponding to $|W_m|$ on $Z_b$.
Then by the same reason for the diagram \eqref{cd3}, we obtain the commutative diagram of meromorphic maps
\begin{equation}\label{cd5}
 \CD
\tilde{Z}_b@>>>Z_b\\
 @V\tilde{\Phi}_m^{G_1}VV @VV{\Phi}_m^{G_1}V\\
\tilde{\mathscr T}_b@>>>\mathscr T_b\\
 \endCD
 \end{equation}
where $\tilde{\mathscr T}_b\to\mathscr T_b$ is the minimal resolution of all singularities of $\mathscr T_b$, and $\tilde{\Phi}_m^{G_1}$ is the meromorphic map  determined by the commutativity of the diagram.
Let $\tilde{\pi}_m:\tilde{\mathscr T}_b\to\Lambda_m$ still denotes the composition of the minimal resolution $\tilde{\mathscr T}_b\to\mathscr T_b$ and $\pi_m:\mathscr T_b\to\Lambda_m$.

Next let $\tilde Z'_b\to\tilde Z_b$ be a sequence of $G_1$-equivariant blowing-ups which eliminates Bs\,$|W_m|$ on $\tilde Z_b$.
We can suppose that the center of the blow-ups are disjoint from $E_1\cup\ol E_1$, where $E_1$ and $\ol E_1$ are the exceptional divisors of the first blow-ups at $C_1$ and $\ol C_1$ as before.
Then we obtain a holomorphic map $\tilde{\Phi}'_b:\tilde Z'_b\to\tilde{\mathscr T}_b$ as the composition  $\tilde{Z}'_b\to\tilde{Z}_b\to\tilde{\mathscr T}_b$.
Again, $\tilde{\Phi}'_b$ is bimeromorphic to $\Phi_m^{G_1}:Z_b\to\mathscr T_b$.
Since our explicit blowing-ups $\tilde Z_b\to Z_b$ are the same as $\tilde Z\to Z$ in \S 3.1, and  the blow-ups $\tilde Z'_b\to \tilde Z_b$ do not touch $E_1\cup\ol{E}_1$, 
the expression \eqref{nb1} for normal bundle is still valid for  $N_{E_1/\tilde{Z}'_b}$.
Hence analogously to \eqref{dis1}, we obtain that the direct image sheaf satisfies
\begin{align}\label{dis2}
\tilde\Phi'_{b*}\mathscr O_{\tilde{Z}'_b}(E_1+\ol{E}_1)
\simeq
N_{E_1/\tilde{Z}'_b}\oplus N_{\ol{E}_1/\tilde{Z}'_b}\oplus\mathscr O_{\tilde{\mathscr T}_b}.
\end{align}
As the meromorphic map associated to the pair $(\tilde{\Phi}'_b, \mathscr O_{\tilde{Z}'_b}(E_1+\ol{E}_1))$ we obtain a meromorphic map
\begin{equation}\label{p2bdle2}
\mu_b: \tilde{Z}'_b\lra\mathbf{P}( N_{E_1/\tilde{Z}'_b}^{\vee}\oplus N_{\ol{E}_1/\tilde Z'_b}^{\vee}\oplus\mathscr O_{\tilde{\mathscr T}_b})
\end{equation}
which is again a  bimeromorphic map onto a conic bundle over $\tilde{\mathscr T}_b$.
We denote this conic bundle by $p_b:X_b\to\tilde{\mathscr T}_b$.
$X_b$ is bimeromorphic to the twistor space $Z_b$.
In particular, $Z_b$ is Moishezon.

%

\begin{prop}\label{prop-d5}
Let $p_b:X_b\to\mathscr{\tilde{T}}_b$ be the above conic bundle which is a projective model of the twistor space $Z_b$ obtained in Theorem \ref{thm-1}. 
Then its discriminant curves consist of the following.
(a) The two distinguished sections $\Gamma$ and $\ol{\Gamma}$ of the morphism $\tilde{\pi}_m:\tilde{\mathscr T}_b\to\Lambda_m$.
(b) The reducible fibers $\tilde{\pi}_m^{-1}(1,\lambda_i,\lambda_i^2,\cdots,\lambda_i^m)$ where $i$ satisfies $r<i<s$ and $l_i>0$.
(c) Irreducible fibers $\tilde{\pi}_m^{-1}(1,\lambda_i,\lambda_i^2,\cdots,\lambda_i^m)$ where $i$ satisfies $r<i<s$ and $l_i=0$.
(d) Irreducible curves belonging to the pull-back of the system $|\mathscr O_{\mathscr T_b}(1)|$ under the minimal resolution $\tilde{\mathscr T}_b\to\mathscr T_b$, where $\mathscr O_{\mathscr T_b}(1)$ is the hyperplane section class with respect to the canonical embedding $\mathscr T_b\subset\mathbf{CP}^{m+2}$.
Further, the number of such discriminant curves  is given by $(n+r-s)$.
\end{prop}

Note that by the assumption \eqref{ineq-key} of Theorem \ref{thm-1}, there actually exists a discriminant curve belonging to (d).
In the proof below, we will show that the conic bundle $p_b:X_b\to\tilde{\mathscr T}_b$ is a deformation of the original conic bundle $p:X\to\tilde{\mathscr T}$, realized inside the $\mathbf{CP}^2$-bundle \eqref{p2bdle2}.
From this viewpoint,
the discriminant curves (a) and  (b) exactly correspond  (a) and (b)  in Proposition \ref{prop-dl1} respectively, while $\{$(c), (d)$\}$ corresponds (c) in Proposition \ref{prop-dl1}.

\bigskip
\noindent {\bf Proof of Proposition \ref{prop-d5}.}
The curves (a), (b) and (c) are discriminant curves by the same reasoning as those for (a), (b) and (c) of Proposition \ref{prop-dl1};
all we have to do is to replace $E_2$ and $\ol{E}_{n+2}$ appeared in the proof (of Proposition \ref{prop-dl1}) by $E_r$ and $\ol{E}_{s}$ respectively.
So we do not repeat the argument here.

To find the discriminant curves in (d) concretely,
fix any integer $i$ satisfying $1<i<r$ or $s<i<n+2$.
The number of such $i$ is $(n+r-s)$.
Then on the twistor space $Z$ of Joyce metric,  the intersection $L_i=S_i^+\cap S_i^-$ was a twistor line contained in $Z^{G_1}$ by Proposition \ref{prop-irf}.
On the present twistor space $Z_b$, although there exist no degree-one divisors $S_i^+$ nor $S_i^-$, the twistor line $L_i$ still makes
as the $G_1$-fixed twistor lines.
Since $L_i$ is not contained in the anticanonical cycle \eqref{ac3} on $Z_b$, the image $\Phi^{G_1}_m(L_i)\subset\mathscr T_b$ makes sense.
In the following we put $\mathscr C_i:=\Phi^{G_1}_m(L_i)$ and show that $\mathscr C_i$ is a hyperplane section of the minitwistor space $\mathscr T_b\subset\mathbf{CP}^{m+2}$.

For this purpose we first show that $\mathscr C_i$ is a curve and that it does not go through the singular point of $\mathscr T_b$.
Since $L_i\subset Z_b^{G_1}$, if $L_i$ intersects the anticanonical cycle \eqref{ac3}, the intersection points are $G_1$-fixed points.
So the intersection point must be on $C_1\cup\ol C_1$ or the singular points of the cycle $C$.
But the latter cannot happen since the twistor line going through the singular point of $C$ is the twistor line $L_j=S_j^+\cap S_j^-$  for $r\le j\le s$ which is different from the $L_i$ (with $1<i<r$ or $s<i<n+2$).
Hence the intersection point of $L_i$ and \eqref{ac3} must be on $C_1$ or $\ol{C}_1$.
But since $C_1$ is a non-isolated $G_1$-fixed locus, any twistor line going through   $C_1$ cannot be contained in the $G_1$-fixed locus (on $Z_b$).
Therefore we obtain that $L_i$ is disjoint from the cycle \eqref{ac3}.
This implies that $L_i\cap$\,Bs\,$|W_m|=\emptyset$.
Hence, since $L_i\cdot mF=2(m-2)>0$,  $\mathscr C_i$ is not a point; namely it is a curve.

Since we have $(\Phi^{G_1}_m)^{-1}(\pi^{-1}_m(\lambda))=\Psi_m^{-1}(\lambda)\in|F|$ for $\lambda\in\Lambda_m$ by the diagram \eqref{cd4} and $L_i\cdot F=2$ on $Z_b$, the curve
$\mathscr C_i$ intersects $\pi^{-1}_m(\lambda)$ transversally at two points for general $\lambda\in\Lambda_m$.
Further, since $L_i\cdot S_j^{\pm}=1$ and $\Phi^{G_1}_m(S_j^{\pm})$ are irreducible components of  a reducible fiber of $\pi_m$ if $l_j>0$, 
$\mathscr C_i$ actually intersects  irreducible components of arbitrary reducible fibers of $\pi_m$.
%
Further, since $L_i$ is away from the cycle \eqref{ac3} and since the conjugate pair of singular points ($P_{\infty}$ and $\ol P_{\infty}$) of $\mathscr T_b$ must be the images of $\sum_{r\le j\le s}C_j$ and $\sum_{r\le j\le s}\ol C_j$, $\mathscr C_i$ does not go through these singularities.
These intersection data directly implies that $\mathscr C_i$ is a hyperplane section of $\mathscr T_b$. 

Finally we show that the conic bundle $p_b:X_b\to\tilde{\mathscr T}_b$ does not have  discriminant curves other than (a)--(d).
For this, recall that the present twistor space $Z_b$ is obtained as a deformation of the twistor space $Z$ of Joyce metric.
Since $W_m=H^0(mF)^G_1$ holds on $Z$ and $Z_b$,
the map $\Phi_m^{G_1}:Z_b\to\mathscr T_b$ is a ($G_1$-equivariant) deformation of 
$\Phi_m^{G_1}:Z\to\mathscr T$.
Further, the explicit  blow-up sequence $\tilde Z_b\to Z_b$ is also a ($G_1$-equivariant) deformation of $\tilde Z\to Z$.
Moreover, by \eqref{bs1}, further blow-ups $\tilde Z'_b\to \tilde Z_b$ can be taken as a $G_1$-equivariant deformation of  $\tilde Z'\to \tilde Z$.
and $Z$,  
the former blow-ups can be taken as a $G_1$-equivariant deformation of the latter blow-ups.
Hence 
the compositions $\tilde Z'_b\to Z_b$ can be supposed to be a $G_1$-equivariant deformation of  $\tilde Z'\to Z$.
Furthermore, since  the bimeromorphic map $\mu_b:\tilde Z_b\to X_b$ and $\mu:Z\to X$ are  
defined as the relative meromorphic maps associated the pair $(\tilde{\Phi}'_b, \mathscr O_{\tilde{Z}'_b}(E_1+\ol{E}_1))$ and $(\tilde{\Phi}', \mathscr O_{\tilde Z'}(E_1+\ol{E}_1))$ respectively and the divisor $E_1+\ol{E}_1$ on $\tilde Z'_b$ corresponds $E_1+\ol{E}_1$ on $Z$ through the deformation, the conic bundle $p_b:X_b\to\tilde{\mathscr T}_b$ is a ($G_1$-equivariant) deformation of the conic bundle $p:X\to\tilde{\mathscr T}$.
Then since $E_r$ and $E_s$ in $Z_b$ correspond to $E_2$ and $E_{n+2}$ in $Z$ respectively, and $S_i^{\pm}$ in $Z_b$ correspond $S_i^{\pm}$ in $Z_b$ for $r<i<s$ under the deformation,  the discriminant curves in (a) and (b) in Proposition \ref{prop-dl1} correspond (a) and (b) of the present proposition.
Also, the discriminant curves in (d)  corresponds a part of (c) (precisely speaking the fibers $\tilde{\pi}_m^{-1}(1,\lambda_i,\cdots,\lambda_i^{m})$ with $1<i<r$ and $s<i<n+2$)
in Proposition \ref{prop-dl1}, and the curves in (c) in the present proposition correspond the remaining curves in (c) of Proposition \ref{prop-dl1}.
Since (a)--(c) in  Proposition \ref{prop-dl1} are all discriminant curves for $p$, it follows that (a)--(d) are also all discriminant curves for $p_b$ as well.
Thus we have shown all the claims of the proposition.
\proofend
\section{Various examples of new Moishezon twistor spaces}
In this section, we shall explain how the results obtained so far produce a number of new Moishezon twistor spaces readily.

\subsection{$U(1)$-equivariant deformations of arbitrary Joyce metrics}
First, we explain a particular (but natural) way for obtaining a pair of a $K$-action on $n\mathbf{CP}^2$ and a $U(1)$-subgroup fixing a sphere, from that on $(n-1)\mathbf{CP}^2$, by means of equivariant connected sum.
For this, we take any effective $K$-action $\rho$ on $(n-1)\mathbf{CP}^2$ and let $S^2_i\subset(n-1)\mathbf{CP}^2$ $(1\le i\le n+1)$ be one of the  $K$-invariant spheres.
Let $K_i\subset K$  be the isotropy $U(1)$-subgroup of $S^2_i$.
Then the $K$-action  has exactly two fixed points on $S^2_i$.
On the other hand we consider a standard $K$-action on $\mathbf{CP}^2$, which is
given by $(z_0,z_1,z_2)\mapsto (z_0,sz_1,tz_2)$, $(s,t)\in K$.
We choose a $U(1)$-subgroup $\{s=1\}$.
The fixed point set of this $U(1)$-action consists of a line $\{z_2=0\}$ and a point $\{(0,0,1)\}$.
Then for any one of the  $K$-fixed  points $p\in S^2_i$ and $q\in \{z_2=0\}$, a  $K$-equivariant connected sum of $(n-1)\mathbf{CP}^2$ and $\mathbf{CP}^2$  at $p$ and $q$ makes sense naturally.
Consequently, we obtain  an effective $K$-action $\rho'$ on $n\mathbf{CP}^2$.
On this $n\mathbf{CP}^2$, there exists a particular $K$-invariant sphere which is  obtained by gluing $S^2_i$ and $\{z_2=0\}$.
Let $K'_i\subset K$ be the isotropy subgroup of the last sphere.
This way, starting from any effective $K$-action $\rho$ on $(n-1)\mathbf{CP}^2$ and any $U(1)$-subgroup fixing a sphere, we naturally obtain an effective $K$-action $
\rho'$ on $n\mathbf{CP}^2$ together with a $U(1)$-subgroup $K'_i$ fixing a sphere.
This operation can be interpreted as $U(1)$-equivariant connected sum 
of $(n-1)\mathbf{CP}^2$ and $\mathbf{CP}^2$ at non-isolated fixed points.
Now we show that all these $U(1)$-actions on $n\mathbf{CP}^2$ have  invariant self-dual metrics which are different from Joyce metrics.
More precisely, we have the following result concerning $U(1)$-equivariant deformations of arbitrary Joyce metrics. 

\begin{thm}\label{thm-3}
Let $n\ge 3$ and consider any effective $K$-action $\rho$ on $(n-1)\mathbf{CP}^2$, and take any one of the $U(1)$-subgroups  $K_i$ ($1\le i\le n+1)$ which fixes one of the $K$-invariant spheres as above.
Then at least one of the following holds.
(i) The $\rho$-invariant Joyce metrics admit a $K_i$-equivariant, non-$K$-equivariant deformation. 
(ii) If $\rho'$ denotes the $K$-action on $n\mathbf{CP}^2$ obtained as a $K$-equivariant connected sum of $(n-1)\mathbf{CP}^2$ and $\mathbf{CP}^2$ at non-isolated fixed points  as above, then  $\rho'$-invariant Joyce metrics admit a  $K'_i$-equivariant, non-$K$-equivariant  deformation.
\end{thm}

\noindent
Proof.
If $\rho$ contains a $U(1)$-subgroup which acts semi-freely (on $(n-1)\mathbf{CP}^2$), then $\rho$-invariant self-dual metrics are LeBrun metrics with torus action \cite{LB93}.
In this case, (i) holds if $n\ge 4$ and (ii) holds if $n=3$.
So in the following we suppose that $\rho$ does not contain such a $U(1)$-subgroup.

Let $Z$ be the twistor space of a $\rho$-invariant Joyce metric on $(n-1)\mathbf{CP}^2$, $S\in |F|^G$ a real smooth member, and $C=\sum_{i=1}^{n+1}(C_i+\ol C_i)$ the $G$-invariant anticanonical curve on $S$ as before.
Then up to a possible permutation for the numbering, the $K$-invariant spheres $S^2_i$ on $(n-1)\mathbf{CP}^2$ are exactly the images of $C_i$ (and $\ol C_i$) under the twistor fibration $Z\to(n-1)\mathbf{CP}^2$.
Now since the $K_i$-action is not semi-free by the above assumption, the two integers $r$ and $s$ (given in Section 4.1) make sense.
When $(n-1)+r-s>0$, (i) holds by Theorem \ref{thm-1}.
When $(n-1)+r-s=0$, as in the proof of Theorem \ref{thm-2}, if $S'\to S$ is a blow-up at a conjugate pair of $G$-fixed points on $C_i\cup\ol C_i$,  
then the twistor space $Z'$ of a Joyce metric on $n\mathbf{CP}^2$ having $S'$ as a member of $|F|^G$ admits a $K'_i$-equivariant, non-$K$-equivariant deformation.
Since a blowing-up is equivalent to taking a connected sum with $\ol{\mathbf{CP}}^2$, by the above choice of $S'\to S$, the $K$-action on $n\mathbf{CP}^2$ is exactly $\rho'$.
This means that $\rho'$-invariant Joyce metrics admit a $K'_i$-equivariant, non-$K$-equivariant deformation.
\proofend

\vspace{2mm}
Let $\delta(n)$ be the number of equivalent classes of $U(1)$-actions on $n\mathbf{CP}^2$ that can be obtained from an effective $K$-action on $n\mathbf{CP}^2$ by taking a $U(1)$-subgroup which fixes one of the $K$-invariant spheres.
(Here, we are considering all effective $K$-actions on $n\mathbf{CP}^2$.)
Then for any inequivalent such $U(1)$-actions on $(n-1)\mathbf{CP}^2$, the resulting $U(1)$-actions on $n\mathbf{CP}^2$ (explained above) are mutually inequivalent.
Therefore, by Theorem \ref{thm-3}, we have the following.
\begin{cor}\label{cor-1}
Let $\delta(n)$ be as above.
Then the number of equivalent classes of $U(1)$-actions on $n\mathbf{CP}^2$ satisfying the following condition is at least $\delta(n-1)$:
(i) the $U(1)$-actions admit invariant self-dual metrics which are different from Joyce metrics,
(ii) their twistor spaces are Moishezon.
\end{cor}


Next we explain how to compute the number $\delta(n)$ for small values of $n$.
In the following we use the notations and conventions we employed in Section 2.
In particular, for a selected subgroup $G_i$ (or $K_i$), $1\le i\le n+2$, we apply a cyclic permutation for the indices of $C_j$'s so that the component $C_i$ becomes $C_1$.
Then we have the sequence $(k_2,k_3,\cdots,k_{n+2})$ of integers
defined by \eqref{f_1}.
The number $k_j$ was the order of the isotropy subgroup of the component $C_j$ of the anticanonical curve.
Namely, the $\mathbf C^*$-subgroup first chosen is acting on the component $C_j$ by weight $k_j$.

First, if $n=0$ (i.\,e. when the 4-manifold is $S^4$), there exists only one effective $K$-action, and the number of $K$-invariant spheres is $2$.
Therefore there are 2 choices of $U(1)$-subgroups which fix the invariant spheres. 
But for both of these 2 subgroups we have $k_2=1$ and they are equivalent actions.
Hence we have $\delta(0)=1$.
Similarly, on $\mathbf{CP}^2$, $K$-action is unique (up to equivalence) and we  have $(k_2,k_3)=(1,1)$ for any of the 3 subgroups.
Hence we again have $\delta(1)=1$.
When $n=2$, it is immediate from the result for the case $n=1$ that the sequence $(k_2,k_3,k_4)$ can take two values $(1,1,1)$ and $(1,2,1).$
Hence we have $\delta(2)=2$.
On $3\mathbf{CP}^2$ the number of $K$-actions is still one, but $(k_2,k_3,k_4,k_5)$ takes 5 values $(1,1,1,1), (1,2,1,1), (1,1,2,1), (1,3,2,1)$ and $(1,2,3,1)$.
The second and third ones, and also the 4-th and 5-th ones, are clearly equivalent, and hence we have $\delta(3)=3$.
 
If $n=4$, the sequence $(k_2,k_3,k_4,k_5,k_6)$ takes 
\begin{align}
(1,1,1,1,1), (1,2,1,1,1), (1,1,2,1,1), (1,1,1,2,1)
\end{align}
which are obtained from $(1,1,1,1)$ by equivariant connected sum (at non-isolated fixed points in general), and 
\begin{align}
(1,1,2,1,1), (1,2,1,1,1), (1,2,3,1,1), (1,2,1,2,1),(1,2,1,1,1)
\end{align}
which are obtained from $(1,2,1,1)$, and
\begin{align}
(1,1,2,3,1), (1,3,2,3,1), (1,2,5,3,1), (1,2,3,4,1),(1,2,3,1,1)
\end{align}
which are obtained from $(1,2,3,1)$. 
Removing equivalent ones, the sequences
\begin{align}
(1,1,1,1,1), (1,2,1,1,1), (1,2,1,2,1),(1,2,3,1,1),(1,3,2,3,1),(1,2,5,3,1),(1,2,3,4,1)
\end{align}
represent all mutually inequivalent $U(1)$-actions.
Hence we obtain $\delta(4)=7$.
By similar computations, we obtain $\delta(5)=15$. 
Although we cannot give an explicit formula for $\delta(n)$, we remark that  it is possible to show that 
\begin{align}\label{quad}
\lim_{n\to\infty}\frac{\delta(n)}{n^2}\ge \frac14.
\end{align}
In particular, $\delta(n)$ increases at least quadratically.
But this  is based on very rough estimate and the actual values seems much bigger.

\subsection{}
In this subsection we apply our result to LeBrun twistor spaces with $K$-action.
For LeBrun's $K$-action on $n\mathbf{CP}^2$ with $n\ge3$, among the $(n+2)$ subgroups, there exist $[(n/2)+2]$ mutually inequivalent $U(1)$-subgroups. 
For these subgroups, the sequence $(k_2,k_3,\cdots,k_{n+2})$  takes the following values:
\begin{align}\label{seq60}
(1,1,\cdots,1),\,(\overbrace{1,2,3,\cdots,n}^n,1),(\overbrace{1,2,3,\cdots,n-1}^{n-1},1,1)
\end{align}
and 
\begin{align}\label{seq61}
 (\overbrace{1,2,3,\cdots,k-1,k}^k,1,\overbrace{n-k,n-k-1,\cdots,2 ,1}^{n-k}),\,\,\left[\frac{n+1}{2}\right]\le k\le n-2.
\end{align}
Here, the sequence $(1,1,\cdots,1)$ corresponds to the semi-free subgroup and in that case, the twistor spaces admit $U(1)$-equivariant deformations which do not preserve $K$-action, provided $n\ge 3$.
Among other sequences, only $(1,2,3,\cdots,n-1,1,1)$ satisfies $n+r-s>0$.
The twistor spaces obtained by equivariant deformations with respect to this $U(1)$-action are exactly the twistor spaces investigated in \cite{Hon07-4}.

The sequences \eqref{seq60} and \eqref{seq61} also mean that from the LeBrun twistor spaces  on $n\mathbf{CP}^2$ with $n\ge 3$, we obtain $[(n/2)+2]$ different minitwistor spaces.
The minitwistor space corresponding to $(1,1,\cdots,1)$ is a smooth quadric in $\mathbf{CP}^3$ (since we have $m=1$), which is of course isomorphic to $\mathbf{CP}^1\times\mathbf{CP}^1$.
The minitwistor spaces arising from the sequences
$(1,2,3,\cdots,n-1,1,1)$ (and 
$(1,2,3,\cdots,n,1)$ also) are the ones studied in \cite{Hon07-4}.
The remaining $[(n/2)-1]$ minitwistor spaces are new (including the author's previous papers).

The $U(1)$-actions generating these new minitwistor spaces satisfy $n+r-s=0$ and hence Theorem \ref{thm-1} does not give $U(1)$-equivariant deformations that do not preserve the $K$-action.
But again by Theorem \ref{thm-3} we obtain a pair of $K$-action on $(n+1)\mathbf{CP}^2$ and a $U(1)$-subgroup for which Theorem \ref{thm-1} gives  $U(1)$-equivariant  but non-$K$-equivariant deformations.
We note that the last $U(1)$-action is not semi-free and the twistor spaces are new.
 
\subsection{}
Next we first display all minitwistor spaces of Joyce metrics that do not have real singularities.
(Recall that by Proposition \ref{prop-sing_T}, the minitwistor spaces always have a conjugate pair of singularities, as long as $m> 1$.)
By Proposition \ref{prop-sing_T}, the minitwistor spaces associated to the $U(1)$-subgroups have real singularity iff some $l_j$ $(1\le j\le n+2)$ satisfies $l_j>1$,where $(l_1,l_2,\cdots,l_{n+2})$ is the sequence of integers defined in Definition \ref{def-m_i}.
So the minitwistor space does not have real singularity only when $l_j=0$ or 1 for all $j$.
It is readily seen that this condition is equivalent to the condition that $k_j=1$ or 2 for all $j$.
Hence the minitwistor spaces associated  to a $U(1)$-subgroups has no real singularities if and only if the isotropy subgroup at any point is either $\{1\}, U(1)$, or $\{\pm1\}$.
The equivalent classes of these $U(1)$ actions are uniquely determined by the number of indices $j$ satisfying $k_j=2$; namely the number of $K$-invariant spheres on which $-1\in U(1)$ acts trivially.
We note that in the sequence $(k_2,k_3,\cdots,k_{n+2})$ a number greater than 1 cannot appear successively, by the effectivity of the actions.
Therefore, on $n\mathbf{CP}^2$, there are exactly $[(n/2)+1]$ kinds of these $U(1)$-actions.
(For example, if $n=7$, these $U(1)$-actions are represented by $(1,1,1,1,1,1,1,1), (1,2,1,1,1,1,1,1), (1,2,1,2,1,1,1,1)$ and $(1,2,1,2,1,2,1,1)$.)

Next we examine equivariant deformations of Joyce metrics with respect to these $U(1)$-actions.
The $U(1)$-action corresponding to $(1,1,\cdots,1)$ is the semi-free $U(1)$-action and therefore they admit $U(1)$-equivalent but non-$K$-equivariant deformations if $n\ge 3$.
For the remaining $U(1)$-actions, $n+r-s$ takes the values
\begin{align}\label{152}
n-2, n-4,n-6,\cdots,2,0
\end{align}
when $n$ is even, and
\begin{align}
n-2, n-4,n-6,\cdots,3,1
\end{align} 
when $n$ is odd.
Therefore, except the final case in \eqref{152}, the twistor spaces admit $U(1)$-equivariant, non-$K$-equivariant deformations by Theorem \ref{thm-1}.
If the number of the $K$-invariant spheres whose isotropy is $\{\pm1\}$   (namely the number of indices satisfying $k_j=2$) is one,
the deformed twistor spaces are exactly the ones investigated in \cite{Hon07-3}.
If the number is greater than one, the  deformed twistor spaces, which are of course Moishezon, are new, to the best of the author's knowledge.
Since the minitwistor spaces of these twistor spaces have no singularity other than the conjugate pair of the quotient singularities, there is a chance that one can find explicit construction of the twistor spaces.

\subsection{}
Finally as another extreme case, for each $n$, we give  minitwistor spaces of Joyce metrics which have a lot of real singularities.
For this, we consider a series of $U(1)$-actions on $n\mathbf{CP}^2$ represented by the following data for $(k_2,k_3,\cdots,k_{n+2})$:

\vspace{2mm}\begin{center}
\begin{tabular}{|c|c|c|c|}\hline
$n$ & $(k_2,k_3,\cdots,k_{n+2})$ & $(l_1,l_2,\cdots,l_{n+2})$ & $m$\\
\hline
2&$(1,2,1)$&$(1,1,1,1)$&2\\
3&$(1,2,3,1)$&$(1,1,1,2,1)$&$3$\\
4&$(1,2,5,3,1)$&$(1,1,3,2,2,1)$&5\\
5&$(1,2,5,8,3,1)$&$(1,1,3,3,5,2,1)$& 8\\
6&$(1,2,5,13,8,3,1)$&$(1,1,3,8,5,5,2,1)$& 13\\
7&$(1,2,5,13,21,8,3,1)$&$(1,1,3,8,8,13,5,2,1)$& 21\\
$\cdots$&$\cdots$&$\cdots$&$\cdots$
\end{tabular}
\end{center}

\vspace{2mm}
\noindent
If $f(n) \,(n=1,2,3,\cdots)$ denotes the Fibonacci sequence so that $f(1)=1,f(2)=1,f(3)=2,f(4)=3,f(5)=5,f(6)=8,f(7)=13,$ etc, then the basic invariant $m$ (defined by Definition \ref{def-m}) for this $U(1)$-action on $n\mathbf{CP}^2$ is given by $f(n+1)$. 
These $U(1)$-actions can be characterized by the property that $m$ attains the maximal value for each $n$, among all $U(1)$-actions on $n\mathbf{CP}^2$ obtained from effective $K$-actions by the restrictions.
If $\mathscr T_n$ denotes the minitwistor space associated to this $U(1)$-action on $n\mathbf{CP}^2$, according to Proposition \ref{prop-sing_T}, $\mathscr T_n$ has real $A_{f(j)-1}$-singularities for all $3\le j\le n$.
So we can say that these minitwistor spaces are the most singular ones among all minitwistor spaces obtained in Section 2.
If $n\ge 4$, these minitwistor spaces are also new, to the best of the author's knowledge.

For these $U(1)$-actions, we always have $n+r-s=0$ and hence Theorem \ref{thm-1} does not generate a $U(1)$-equivariant deformation which does not preserve $K$-action.
But as we have frequently done, by taking the $K$-equivariant connected sum with $\mathbf{CP}^2$ at the fixed sphere, the Joyce metrics on $(n+1)\mathbf{CP}^2$ invariant under the resulting $K$-action admit a $U(1)$-equivariant deformation that does not preserves the $K$-action.
The corresponding twistor spaces is Moishezon twistor space with $\mathbf C^*$-action whose quotient space is $\mathscr T_n$ by Theorem \ref{thm-2}.
These twistor spaces, whose detailed structure can be also derived from Proposition \ref{prop-d5}, are new in view of the $U(1)$-actions on $n\mathbf{CP}^2$.
However, contrary to the twistor spaces in \S 5.3, it seems difficult to obtain an explicit construction of these twistor spaces, because the projective models have a lot of complicated singularities.

\section{Appendix} 
In  Section 2 we obtained minitwistor spaces of  Joyce metrics as an image of the meromorphic map associated to some explicit linear system on the twistor spaces, and showed that they can be regarded as quotient spaces of the twistor spaces by $\mathbf C^*$-action.
In this appendix we show that these minitwistor spaces are actually `canonical' quotient spaces. 

For this we first recall a result of  Fujiki \cite{F78-2} concerning quotient spaces under holomorphic actions of  a complex Lie group on compact complex manifolds.
For simplicity we explain in simple situation which is enough for our purpose.
Let $X$ be a compact complex manifold in Fujiki class $\mathscr C$  and suppose that a Lie group $\mathbf C^*$ is acting holomorphically on $X$.
Further, suppose that the action has at least one fixed point.
Let $D_X$ be the Douady space of $X$, which is a complex space parametrizing all complex subvarieties in $X$.
By the assumption $X\in\mathscr C$, all irreducible components of $D_X$ are compact \cite{F78-1}.
The $\mathbf C^*$-action on $X$ naturally induces a holomorphic $\mathbf C^*$-action on $D_X$.
The fixed point set $D_X^{\mathbf C^*}$ of the last action parametrizes all $\mathbf C^*$-invariant subvarieties in $X$.

\begin{prop}\label{prop-fujiki}(Fujiki \cite[Lemma 4.2]{F78-2})
There exists a unique irreducible component $Y$ in $D_X^{\mathbf C^*}$ satisfying the following property:
There exists a dense Zariski open subset $U\subset Y$ such that for any $u\in U$, the corresponding $\mathbf C^*$-invariant subvariety in $X$ is the closure of an orbit of the $\mathbf C^*$-action.
\end{prop}

If  $(X\times D_X\supset)\,\mathscr W\to D_X$ denotes the universal family and $\mathscr W_Y\to Y$ its restriction to the subspace $Y$ in the proposition, the restriction of the natural projection $\mathscr W\to X$ to $\mathscr W_Y$ is bimeromorphic. Hence by composition with $\mathscr W_Y\to Y$, we obtain a meromorphic map $f:X\to Y$.
By construction, for any $u\in U\subset Y$, $f^{-1}(u)$ is the closure of an orbit of the $\mathbf C^*$-action.
Since $D_X$ and $D_X^{\mathbf C^*}$ are canonically determined from the space $X$ and the $\mathbf C^*$-action, we call the uniquely determined space $Y$ in the proposition the {\em canonical quotient space} (by the $\mathbf C^*$-action on $X$).

We go back to the situation we have been considering.
Let $Z$ be the twistor space of a Joyce metric on $n\mathbf{CP}^2$ and $G_1$ ($\simeq\mathbf C^*$) a subgroup of $G$ ($\simeq\mathbf C^*\times\mathbf C^*$) fixing (any one of) a component $C_1$ of the cycle $C$.
Let $\mathscr T$ be the minitwistor space of the Joyce metric with respect to $G_1$ (in the sense of Definition \ref{def-mt}).
We will show the following.

\begin{prop}\label{prop-cq}
The minitwistor space $\mathscr T$ is isomorphic to the canonical quotient space by the $G_1$-action on $Z$.
\end{prop}

For the proof, we need the following

\begin{lemma}\label{lemma-cq}
As before let $\Phi_m^{G_1}:Z\to\mathscr T$ be the meromorphic map associated to the linear system $|W_m|$ (see Definition \ref{def-w}).  Then we have the following.
(i) There is a Zariski open subset $U\subset\mathscr T$ such that for any $u\in U$ the fiber $(\Phi_m^{G_1})^{-1}(u)$ is the closure of an orbit of the $G_1$-action.
(ii)  $\Phi_m^{G_1}$ does not contract any divisor to a point.
\end{lemma}

\noindent Proof.
For (i) recall that we have the commutative diagram \eqref{cd2}.
We set $\Lambda_m^{\circ}:= \Lambda_m\backslash\{\lambda_1,\cdots,\lambda_{n+2}\}$ and  $U:=\pi_m^{-1}(\Lambda_m^{\circ})\backslash\{p_{\infty},\ol p_{\infty}\}$.
If $\lambda\in \Lambda^{\circ}$,  $S_{\lambda}:=\Psi_m^{-1}(\lambda)$ is an irreducible and smooth member of the pencil $|F|^G$.
Further by Propositions \ref{prop-mt2} and \ref{prop-mt3} (ii), the restriction $\Phi_m^{G_1}|_{S_{\lambda}}$ is identical to the  quotient morphism \eqref{quot_S}. 
Since the map \eqref{quot_S} has exactly two singular (i.e.\,reducible) fibers, and these are mapped to $p_{\infty}$ and $\ol{p}_{\infty}$, we obtain that $U$ satisfies  the claim (i).

For (ii), again by the commutative diagram \eqref{cd2}, it suffices to show that $\Phi_m^{G_1}$ does not contract any irreducible components of a member in  $|F|^G$ to a point.
As in the above proof for the claim (i), this is obvious for irreducible members.
So it remains to see that $\Phi_m^{G_1}(S_i^+)$ and $\Phi_m^{G_1}(S_i^-)$ are not points for any $1\le i\le n+2$, where $S_i^++S_i^-$ are reducible members of $|F|^G$.
For this, we note that the meromorphic map $\Phi_m^{G_1}|_{S_i^+}$ is exactly the rational map associated to a linear system $|W_m|_{S_i^+}|$, where $W_m|_{S_i^+}$ is the image of the subspace $W_m\subset H^0(Z,mF)$ under the restriction map $H^0(Z,mF)\to H^0(S_i^+,mF|_{S_i^+})$.
By the definition of $W_m$, $|W_m|$ is generated by $Y, \ol Y$ and $|V_m|$.
Hence $|W_m|_{S_i^+}|$ is generated by $Y|_{S_i^+},\ol Y|_{S_i^+}$ and $mS|_{S_i^+}$, where $S$ is a smooth member of $|F|^G$.
Then by Proposition \ref{prop-mt3} (iv), at least one of $S_i^+\not\subset Y$ or $S_i^+\not\subset \ol Y$ holds. 
Further, of course $S_i^+\not\subset S$ holds. 
Hence the linear system $|W_m|_{S_i^+}|$ has at least two different members.
Therefore $\Phi_m^{G_1}(S_i^+)$ is not a point.
Then by reality, $\Phi_m^{G_1}(S_i^-)$ is not a point too.
\proofend

\vspace{2mm}\noindent
Proof of Proposition \ref{prop-cq}.
Let $\alpha:\ol Z\to Z$ be a sequence of $G_1$-equivariant blowing-up which eliminates the indeterminacy of $\Phi_m^{G_1}:Z\to\mathscr T\,$ so that the composition $\ol Z\to Z\to\mathscr T$ is a $G_1$-equivariant morphism.
Then by Hironaka's flattening theorem \cite[Corollary 1]{H75}, there is a sequence of blowing-ups $\mathscr T'\to\mathscr T$, for which if $\Phi':Z'\to\mathscr T'$ denotes the strict transform of the fiber product $\ol Z\times_{\mathscr T}\mathscr T'\to\mathscr T'$ and $\beta:Z'\to\ol Z$ denotes the natural projection, the morphism $\Phi'$ becomes flat:
\begin{equation}\label{cd6}
 \CD
Z'@>{\beta}>>\ol Z@>{\alpha}>>Z\\
 @V{\Phi'}VV @VVV @VV{\Phi}_m^{G_1}V\\
{\mathscr T}'@>>>\mathscr T@=\mathscr T\\
 \endCD
 \end{equation}By the universality of the Douady space and the flatness of $\Phi'$, $\mathscr T'$ can be regarded as a complex subvariety of the Douady space $D_{Z'}$.
Further, by Lemma \ref{lemma-cq} (i), this subvariety $\mathscr T'$ must be the unique component of $D_{Z'}^{G_1}$ given in Proposition \ref{prop-fujiki}.

We write $\gamma=\alpha\circ\beta$ and let $\gamma_*:D_{Z'}\to D_Z$ be the holomorphic map induced by $\gamma$.
Since  the left half of the diagram \eqref{cd6} automatically becomes $G_1$-equivariant,
$\gamma_*$ is obviously $G_1$-equivariant.
Further, since general fibers of $Z'\to\mathscr T'$ is mapped to a general fiber of $\Phi_m^{G_1}:Z\to\mathscr T$, the image $\gamma_*(\mathscr T')$ must equal to the unique component of $D_Z^{G_1}$ in Proposition \ref{prop-fujiki}.
Namely $\gamma_*(\mathscr T')$ coincides with the canonical quotient space of $Z$ by the $G_1$-action.
Take any $y\in \mathscr T$ and let $y'_1$ and $y'_2$ be points on $\mathscr T'$ which are mapped to the same point $y$ under $\mathscr T'\to\mathscr T$. 
Then by the commutativity of the diagram \eqref{cd6}, $\gamma((\Phi')^{-1}(y_1))$ and $\gamma((\Phi')^{-1}(y_2))$ 
are contained in $(\Phi_m^{G_1})^{-1}(y)$.
Now since $(\Phi_m^{G_1})^{-1}(y)$ is a curve by Lemma \ref{lemma-cq} (ii), we have $\gamma((\Phi')^{-1}(y_1))=\gamma((\Phi')^{-1}(y_2))=(\Phi_m^{G_1})^{-1}(y)$. This means $\gamma_*(\mathscr T')\simeq\mathscr T$, as required.
\proofend

\small
\vspace{13mm}
\hspace{7.5cm}
$\begin{array}{l}
\mbox{Department of Mathematics}\\
\mbox{Graduate School of Science and Engineering}\\
\mbox{Tokyo Institute of Technology}\\
\mbox{2-12-1, O-okayama, Meguro, 152-8551, JAPAN}\\
\mbox{{\tt {honda@math.titech.ac.jp}}}
\end{array}$

\end{document}